\documentclass[10 pt]{article}

\newtheorem{theorem}{Theorem}[section]

\newtheorem{proposition}[theorem]{Proposition}

\newtheorem{lemma}[theorem]{Lemma}
\newtheorem{definition}{Definition}

\def \dis {\displaystyle}

\newcommand{\myfrac}[2]{{\displaystyle\frac{#1}{#2}}}
\newcommand{\dpart}[2]{{\myfrac{\partial #1}{\partial #2}}}
\newcommand{\dt}[1]{{\myfrac{\partial #1}{\partial t}}}

\newcommand{\arglim}[2]{{\displaystyle{\mathop{\scriptstyle{#1}}_{#2}}}}

\newcommand{\R}{I\!\!R}
\newcommand{\N}{I\!\!N}

\newcommand{\eps}{{\varepsilon}}

\newcommand{\B}{\mathcal{B}}
\newcommand{\C}{\mathcal{C}}
\newcommand{\F}{\mathcal{F}}
\newcommand{\K}{\mathcal{K}}
\newcommand{\G}{\mathcal{G}}
\newcommand{\X}{\mathcal{X}}

\newcommand{\iGG}{{\stackrel{\scriptscriptstyle\rm o}{\G}}}

\newcommand{\GO}{\mathcal{O}}
\newcommand{\oGO}{{\overline{\GO}}}
\newcommand{\dGO}{{\partial \GO}}

\newcommand{\ue}{{u_\eps}}

\def \cx {\hat{x}}
\def \cy {\hat{y}}
\def \ct  {\hat{t}}
\def \cs  {\hat{s}}

\begin{document}

\title{Large Deviations Principle by viscosity solutions: the case of diffusions with oblique Lipschitz reflections.}

\author{Magdalena Kobylanski\footnote{Universit\'e Paris-Est; magdalena.kobylanski@univ-mlv.fr}}

\maketitle

\begin{abstract}
We establish a Large Deviations Principle for stochastic processes with Lipschitz continuous oblique reflections on regular domains. The rate functional is given as the value function of a control problem and is proved to be good. The proof is based on an original viscosity solution approach. The idea consists in interpreting the probabilities as the solutions of some PDEs, make the logarithmic transform, pass to the limit, and then identify the action functional  as the  solution of the limiting equation. 
\end{abstract}

\vspace{13mm}


\vspace{3mm}



\section{Introduction}

According to the terminology of Varadhan \cite{V84}, a sequence $(X^\eps)$ of random variables with values in a metric space $({\mathcal X}, d)$ satisfies a {\em Large Deviations Principle} (LDP in short) if 

\vspace{2mm}

\noindent{\em
There exists a lower semi-continuous functional $\lambda: \X \to [0,\infty]$ such that for each Borel measurable set $\G$ of $\X$
\begin{enumerate}
\item[\rm{(I)}] $\displaystyle \limsup_{\eps\rightarrow 0} \left\{-\eps^2 \ln {\bf
P}[X^{t,x,\eps}\in
\G]\right\}\leq \inf_{g\in\stackrel{\scriptscriptstyle\rm o}{\G}}\lambda_{t,x} (g)$  (LDP's upper bound),
\item[\rm{(II)}] $\displaystyle
\inf_{g\in\overline\G}\lambda_{t,x} (g)\leq\liminf 
_{\eps\rightarrow 0}
\left\{-\eps^2 \ln {\bf
P}[X^{t,x,\eps}\in
\G]\right\}$  (LDP's lower bound),
\end{enumerate}
$\lambda$ is called the {\em rate functional} for the large Deviations principle (LDP). A rate functional is {\em good} if for any $a\in[0,\infty)$, the set $\{g\in \X: \lambda(g)\le a\}$ is compact.}\\

We refer the reder to the books \cite{A80}, \cite{DP83},  \cite{FreidlinWentzell84}, \cite{S84},\cite{V84}, \cite{DZ}, \cite{DE} for the general theory, references and different approches to Large Deviations. 

Partial Differential Equations (in short PDEs) methods have been applied to establish different types of Large Deviations estimates starting from \cite{Flem78}. The idea consists in interpreting the probabilities as the solutions of some PDEs, make the logarithmic transform, pass to the limit, and then identify the action functional  as the  solution of the limiting equation. The 
notion of viscosity solutions (cf. \cite{CL83},  \cite{PL2-83}-\cite{PL2-85}, \cite{Crandall-Ishii-Lions92}) appeared to be particularly adapted to this problem. Indeed,  the half-relaxed semi-limit method  (cf. \cite{BP2}) allows to pass to the limit very easily, moreover the notion of strong uniqueness for viscosity solution allows to identify the solution of the limiting equation with the action functional. A number of Large 
Deviations results have been proved by using this method  \cite{EvansIshii85}, \cite{FlemSouga86}, \cite{BP2}, \cite{BP3}, \cite{BB97}, \cite{AD99},  \cite{Ph03}. However it was a long lasting critic to this method not to provide the  general Large Deviations Principle. The aim of this work was to overcome this gap. We carry out this method in order to establish a LDP for small diffusions with oblique Lipischitz continuous direction of reflections which explains the technicity. This result is new to the best of our knowledge. Our method which was developed in \cite{K98}, seems very efficient and we hope it gives a new insight.

 Recently \cite{FK06} came to our knowledge. This book shows also, in a very general setting, that viscosity solutions are an adapted tool in order to establish LDPs. 

\vline

To be more specific we introduce the precise mathematical formulation of the problem.

\vline 

Let $\GO$ be a smooth open bounded subset of $\R^d$. For $(t,x)\in \R^+\times \oGO$, we consider the oblique reflection problem

\begin{equation}
\label{Eo}
\left\{\begin{array}{l}
dX_s=b(s,X_s)ds-dk_s,\quad X_s\in\overline \GO\
(\forall s> t), \vspace{2mm}\\
{\displaystyle k_s=\int_t^s{\bf 1}_{\partial
\GO}(X_\tau)\gamma(X_\tau)d|k|_\tau\ (\forall s>
t), \quad X_t=x.}
\end{array}\right.
\end{equation}
where $b$ is a continuous $\R^d$-valued function defined on $\R^+\times \oGO$ and $\gamma$ is a $\R^d$-vector field defined on $\dGO$. The solutions of problem (\ref{Eo}) are pairs $(X,k)$ of continuous functions from $[t,\infty)$ to $\oGO$ and $\R^d$ respectively such that $k$ has bounded variations, and $|k|$ denotes the total variation of $k$.

We shall denote by $n(x)$  the unit outward normal to $\dGO $  at $x$,  and assume that
\begin{equation}\label{gamma}
\left\{\begin{array}{l} 
\gamma:\R^d\to \R^d\mbox{ is a Lipschitz continuous function}  \vspace{2mm}\\ \quad
\mbox{and } \exists c_0  \; \forall  x\in\dGO,\; \gamma(x)\cdot n(x)\geq c_0>0.
\end{array}\right.
\end{equation}

When $b$ is Lipschitz continuous, $\gamma$ satisfies condition (\ref{gamma}) and $\GO$ is smooth, the existence of the solutions of  $(\ref{Eo})$ is given as a particular case of the results of  Lions and Sznitman \cite{LS84}  and the uniqueness is a corollary   of the result of Barles and Lions \cite{BL}. For more general domains existence and 
uniqueness of solutions  of $(\ref{Eo})$  is given as a particular case of Dupuis and Ishii \cite{DI93}. The reader can also use the results given in Appendix B.

\vline

 Let 
$(\Omega,\F,(\F_t)_{t\geq 0},{\bf P})$
be a probability space which satisfies the {\em usual conditions} and $(W_t)_{t\geq 0}$ be a standard Brownian
motion with values in $\R^m$.
Consider for each $\varepsilon>0$, $t\geq 0$,
$x\in \oGO$, the following stochastic differential equation

\begin{equation}
\label{Eoe}
\left\{\begin{array}{l}
dX_s^\eps=b_\eps(s,X_s^\eps)ds+\eps
\sigma_\eps(s,X_s^\eps)dW_s-dk_s^\eps,\quad X_s^\eps\in\overline \GO\
(\forall s> t), \vspace{2mm}\\
{\displaystyle k_s^\eps=\int_t^s{\bf 1}_{\partial
\GO}(X_\tau^\eps)\gamma(X_\tau^\eps)d|k^\eps|_\tau\ (\forall s>
t), \quad X_t^\eps=x,}
\end{array}\right.
\end{equation}
where $\sigma$ is continuous $\R^{d\times m}$-valued. A strong solution of (\ref{Eoe})    is a couple $(X^\eps_s,k^\eps_s)_{s\geq t}$ of $(\F_s)_{s\geq t}$-adapted 
processes which have almost surely continuous paths and such that 
$(k^\eps_s)_{s\geq t}$ has
almost surely bounded variations, and $|k^\eps|$ denotes its total variation.

\vline

 Let us now make some comments about this reflection problem. Consider equation (\ref{Eoe}) in the case when $\varepsilon=1$.

This type of stochastic differential equations has been solved by using the Skorokhod map by Lions and Sznitman  in  \cite{LS84} in the case when $\GO$ belongs to a very large class of admissible open subsets and the direction of reflection is the normal direction $n$, or when $\GO$ is smooth and $\gamma$ is of class $C^2$. This problem was 
also deeply studied by Dupuis and Ishii \cite{DI91}, \cite{DI91b}, \cite{DI93}. When $\GO$  is convex these authors proved in \cite{DI91} that the Skorokhod map is Lipschitz continuous  even when trajectories may have jumps. As a corollary, this result gives existence and uniqueness of the solution of the stochastic equation (\ref{Eoe}) and  provides the Large Deviations estimates as well. Dupuis and Ishii also proved in \cite{DI93} the existence of the solution of equation  (\ref{Eoe}) in the following cases: either $\gamma$ is $C^2$ and $\GO$ has only an exterior cone condition, or $\GO$ is a finite intersection of $C^1$ regular bounded domains $\GO_i$ and $\gamma$   is Lipschitz continuous at points $x\in \dGO$ when $x$ belongs to only one $\dGO_i$  but  when $x$ is a corner point, $\gamma(x)$ can even be 
multivaluated. A key ingredient is the use of test functions that Dupuis and Ishii build in \cite{DI90}, \cite{I91} and \cite{DI91b} in order to study oblique derivative problems for fully nonlinear second-order elliptic PDEs on nonsmooth domains.

Let us point out that these type of diffusions with oblique reflection in domains with corners arise as rescaled  queueing networks and related systems with feedback (see \cite{AD99} and the references within).

\vline 

We study in the present paper Large Deviations of (\ref{Eo}) under the simpler condition of a domain without corners. More precisely we suppose that
\begin{equation}\label{Dom}
\GO \mbox{ is a }  W^{2,\infty}\mbox{ open bounded set of } \R^d.
\end{equation}

\vline

Let us precise now what is the regularity we require on the coefficients $b$,  $\sigma$ and  $b_\eps$, $\sigma_\eps$ and how $b_\eps$  and $\sigma_\eps$ are supposed to converge to $b$ and $\sigma$.

For all $\eps>0$, let $ b_{\eps}, b\in {\mathcal C}([0,+\infty)\times 
\oGO;\R^d),$ $
\sigma_{\eps}, \sigma\in
 {\mathcal C}([0,+\infty)\times\oGO;\R^{d\times m})$. And assume that  for each $T>0$, there exists a constant $C_T$ such that for all $\eps>0$,  for all $t$ $\in $ $[0,T]$, for all $x,x'$ $\in $ $\oGO$ one has 
\begin{equation}\label{bsigma1}
\begin{array}{c}|b(t,x)-b(t,x')|,\;|b_\eps(t,x)-b_\eps (t,x')|\le C_T|x-x'|, \\
\|\sigma(t,x)-\sigma(t,x')\|,\;\|\sigma_\eps(t,x)-\sigma_\eps (t,x')\|\le C_T|x-x'|.
\end{array}
\end{equation} 

  We also assume that  
\begin{equation}\label{bsigma2}
(b_\eps), (\sigma_\eps) \mbox{ converge uniformly to } b \mbox{ and } \sigma \mbox{ on } 
[0,T]\times \oGO.
\end{equation}

By \cite{DI93} for all $\eps>0$ and for
all $(t,x)\in[0,T]\times\oGO$, there
exists a unique solution $(X^{t,x,\eps},k^{t,x,\eps})$ of
(\ref{Eoe}) on $[t,T]$. Morever $X^{t,x,\eps}$ converges in probability to the solution $X^{t,x}$ of (\ref{Eo}) when $\eps$ converges to $0$. Obtaining the Large Deviations estimates provides the rate of this convergence.

\vline 

We now turn to the definition of the rate functional $\lambda$. It is defined under conditions (\ref{gamma})-(\ref{Dom})-(\ref{bsigma1}) as  the value function of a non standard contole problem of a deterministic 
differential equation with 
$L^2$ coefficients and with oblique reflections.

More precisely,
let $(t,x)\in [0,T]\times\oGO$ and $\alpha\in L^2(t,T;\R^m) $ and consider equation
\begin{equation}
\label{Y}
\left\{\begin{array}{l}
dY_s=\left( b(s,Y_s)-\sigma(s,Y_s)\alpha_s\right)ds-dz_s, \quad Y_s\in\oGO,(\forall s> t), \vspace{2mm}\\
z_s= \displaystyle   \int_t^s {\bf 
1}_{\dGO}(Y_\tau)\gamma(Y_\tau)d|z|_\tau,(\forall s> t), \quad Y_t=x.
\end{array}\right.
\end{equation}
We prove in  Appendix B that there exists a unique solution $(Y^{t,x,\alpha}_s,z^{t,x,\alpha}_s)_{s\in[ t,T]}$ of (\ref{Y}), and we study the regularity of $Y$ with respect to $t,x,\alpha$ and $s$.

\vline 

In the following we note 
$$\X= \C([0,T];\oGO) \quad \mbox{  and for } g\in \X, \quad\|g\|_{\X}=\sup_{t\in[0,T]} |g(x)|.$$

We  make the following abuse of notations. For $\G\subset$ $\X$ and for $g\in$ $\C([t,T];\oGO)$ for some $t\in[0,T]$, we write $g\in \G$ if there exists a function in $\G$ whose restriction to $[t,T]$ coincides with $g$. 

\vline

 For all $g\in \X$,  we define $\lambda_{t,x}(g)$ by 
\begin{equation}
\label{lambda}
\lambda_{t,x}(g)=\inf \left\{ \frac{1}{2}\int_t^T|\alpha_s|^2ds;\ 
\alpha\in L^2(t,T),\;  Y^{t,x,\alpha} = g \right\}.
\end{equation}
Note that $\lambda_{t,x}(g)\in[0,+\infty]$. 

\vline 

The main result of our paper is the proof of
the full Large Deviations type estimates for (\ref{Eoe}), as well as the identification of the rate functional which is proved to be  good.

\begin{theorem} \label{T}
Assume (\ref{gamma})-(\ref{Dom})-(\ref{bsigma1})-(\ref{bsigma2}). For each $(t,x)\in[0,T]\times \oGO$, and  $\varepsilon>0$,  denote
by $X^{t,x,\eps}$ the unique solution of (\ref{Eoe}) on
$[t,T]$. Consider $\lambda_{t,x}$ defined by (\ref{lambda}). Then 
$(X^{t,x,\eps})_\eps$ satisfies a {\em Large Deviations Principle} with rate functional $\lambda_{t,x}$. Moreover the rate functional is good.
\end{theorem}

As far as the partial differential equations are concerned, we use
the notion of viscosity solutions. We shall not recall the
classical results of the theory of viscosity solutions here
and we refer the reader to
M.G.~Crandall, H.~Ishii and P.-L.~Lions~\cite{Crandall-Ishii-Lions92} (Section 7 for 
viscosity
solutions of second order Hamilton-Jacobi equations), to W.H.
Fleming and H.M. Soner \cite{FlemSoner93} (Chapter 5 for stochastic controlled 
processes) and to G.~Barles~\cite{Barles94} (Chapter 4 for
viscosity solutions of first order Hamilton-Jacobi equations with
Neumann type boundary conditions and Chapter 5 Section 2 for
deterministic controlled processes with reflections).

\vline

The paper is organised as follows. In section 2, we prove  first that assertion (I)  amounts to the proof of this upper bound for a ball $\B$ (assertion (A1)). Second,  we prove that if the rate is good,  assertion (II) amounts   to prove the lower bound for a finit intersection of complementaries of balls (assertion (A2)). Finally, we prove that the fact that the rate is good holds true if  a stability result  holds true for equation (\ref{Y}) (assertion (A3)). In section 3, we give the proof of (A1), and we finish in section 4 by the proof of (A2). 

An important Appendix follows. It includes, in Appendix B, the study of equation (\ref{Y}) and the proof of (A3). In Appendix C, we  study  different mixed optimal control-optimal single or multiple stopping times problems and  we caracterize particular value functions as the minimal (resp. maximal) viscosity supersolution (resp. subsolution) of the related obstacle problems. These caracterizations are important in order to establish (A1) and (A2). Eventually, in Appendix D,  we prove a strong comparison result for viscosity solutions of an obstacle problem with Neumann boundary conditions and quadratic growth in the gradient in the case of a continuous obstacle. This result is needed in the proof of the caracterization of the value functions mentioned above. This long and technical Appendix begins in Appendix A, by the construction of an appropriate test function which is usefull in order to establish the results concerning equation (\ref{Y}) (Appendix B) and the uniqueness result (Appendix D).

\section{A preliminary result}\label{LDE}

We now define the action functional.  For each $(t,x)$ $\in$ $[0,T]$ $\times$ $\oGO$, and for each $\G$ $\subset$ $\X$ let us define $\Lambda_{t,x}(\G)$ as follows:
\begin{equation}\label{Lambda}
\Lambda_{t,x}(\G)= \inf\left\{ \frac{1}{2}\int_t^T|\alpha_s|^2ds;\ 
\alpha\in L^2(t,T),  Y^{t,x,\alpha}\in \G  \right\},
\end{equation}
where $Y^{t,x,\alpha}$ is defined by (\ref{Y}).

It is straightforward that $\Lambda_{t,x}$ is decreasing along increasing sequences of sets and that 
$$ \Lambda_{t,x}(\G)= \inf_{g\in \G}\lambda_{t,x}(g) \quad \mbox{ and } \quad \lambda_{t,x}(g)=\Lambda_{t,x}(\{g\}). $$
We use the following notation: for  $g_0\in \X$ and $r>0$ we denote by $\B(g_0,r)$ the ball of center $g_0$ and of radius $r$ that is $\B(g_0,r)=\{g\in \X, \; \|g-g_0\|_\infty <r\}$.

\vspace{3mm}

We consider the following assertions.

\begin{enumerate}
\item[(A1)]
for all $ g\in \X,$ and $r>0,$
$$\limsup_{\varepsilon \to 0}\left\{ -\varepsilon^2 \ln P[X^{t,x,\eps}\in\B(g,r)]\right\}\le  \Lambda_{t,x}(\B(g,r)), $$
\item[(A2)] 
for all  $i\in \N$,  $  g_i\in \X$, and  $r_i>0$, setting $\B_i=\B(g_i,r_i)$, one has for all $ N\in \N^*,$
$$\quad \liminf_{\eps\to 0}\left\{-\eps^2\ln  P\left [X^{t,x,\eps}\in \bigcap_{i=1}^N
\B_i^c \right]\right\}\ge \Lambda_{t,x}\left( \bigcap_{i=1}^N
\B_i^c\right),
$$
\item[(A3)]
for all $\alpha_n, \alpha\in L^2(0,T;\R^m)$, \\
 if $\alpha_n \rightharpoonup \alpha $ weakly in $L^2$ then  $ \|  Y^{t,x,\alpha_n}- Y^{t,x,\alpha}\|_{\X}\to 0.$\\
\end{enumerate}

The following proposition shows that Theorem \ref{T} reduces to assertions (A1), (A2) and (A3).
\begin{proposition}\label{Pred} For all $(t,x)\in[0,T]\times\oGO$ one has, 

$\quad  \;$(i)$\quad  \;$(A1) implies  (I),

$\quad  \,$(ii)$\quad  \,$If the rate is good then (A2)  implies   (II).

$\quad$(iii)$\quad$   (A3)  implies  that  the rate functional $\lambda_{t,x}$ defined by (\ref{lambda}) is good.
\\
In particular, the proof of Theorem \ref{T} amounts to the proofs of (A1), (A2), (A3).
\end{proposition}

\noindent{\em Proof: } First let us prove (i).

Consider a measurable set $\G\subset \X$.
For all $g\in\iGG$, there exists $r>0$ such that
$\B(g,r)\subset\iGG$. Then
$$-\eps^2\ln P[X^{t,x,\eps}\in\G]\leq 
-\eps^2\ln P[X^{t,x,\eps}\in\iGG]\leq
-\eps^2\ln P[X^{t,x,\eps}\in \B(g,r)].$$
Hence, by (A1) we obtain 
$$\limsup_{\eps\to 0}\left\{-\eps^2\ln  P[X^{t,x,\eps}\in\G]\right\}\le \limsup_{\eps\to 0}\left\{-\eps^2\ln  P[X^{t,x,\eps}\in\B(g,r)]\right\}$$
$$ \le \Lambda_{t,x}(\B(g,r))\leq \lambda_{t,x}(g),$$
and we conclude by taking the infimum over all $g\in\iGG$. This completes the proof of point (i).

\vline

\noindent Let us prove (ii).

Fix
$a<\Lambda_{t,x}(\overline\G)$ and put 
$\K=\{g\in \C([t,T];\oGO),\; \lambda_{t,x}(g)\leq a\}$. Note that $\Lambda_{t,x}(\K^c)\geq a$.
One also clearly has $\K\subset
\overline\G^c$, indeed if $g\in \K$ one has $\lambda_{t,x}(g)\le a$ and if $g\in\overline\G$ then $\lambda_{t,x}(g)>a$ hence $\K\cap \oGO=\emptyset$. As a concequence, for all $g\in\K$, there exists $r>0$ such that
$\B(g,r)\subset \overline\G^c$. Since, by (T2), $\K$ is compact,
there exists a finite number $N$ of
$\B_i=\B(g_i,r_i)$ with $g_i\in \K$ and
$r_i>0$,  such that
$\displaystyle \K\subset \bigcup_{i=1}^N\B_i\subset \overline\G^c.$
 In view of the fact that $\Lambda_{t,x}$ is decreasing
along increasing sequences of sets, this yields, passing to the complementaries that
$$   \Lambda_{t,x}\left( \bigcap_{i=1}^N
\B_i^c\right)\ge \Lambda_{t,x}(\K^c)\ge a.$$
It is clear that 
$$\displaystyle \liminf_{\eps\to 0}\left\{-\eps^2\ln P\left [X^{t,x,\eps}\in 
\G\right ]\right\} \ge \liminf_{\eps\to 0}\left\{-\eps^2\ln P\left [X^{t,x,\eps}\in 
\bigcap_{i=1}^N
\B_i^c\right ]\right\}. $$
But, by (A3) we also have 
$$\liminf_{\eps\to 0}\left\{-\eps^2\ln P\left [X^{t,x,\eps}\in 
\bigcap_{i=1}^N
\B_i^c\right ]\right\}\geq  \Lambda_{t,x}\left( \bigcap_{i=1}^N
\B_i^c\right),  $$
hence we have shown that for all $a<\Lambda_{t,x}(\overline\G)$,
we have 
$$ \liminf_{\eps\to 0}\left\{-\eps^2\ln P\left [X^{t,x,\eps}\in 
\G\right ]\right\} \ge a.$$
Passing to the limit when $a$ tends to $ \Lambda_{t,x}(\overline\G)$, we have completed  the proof of (ii).

\vline

\noindent Let us prove (iii).  

We suppose that the rate functional $\lambda_{t,x}$  defined by (\ref{lambda}) satisfies (A3).
Fix $(t,x)\in [0,T]\times\oGO$, and $a\in\R$. Put $\K=\{g\in \X,\;\lambda_{t,x}(g)\leq a\}$.
Let $(g_n)_{n\in\N}$ be a sequence of $\K$. Then, for all
$n$, there exists $\alpha_n\in  L^2(t,T)$ such that
$Y^{t,x,\alpha_n}=g_n$ and
$$\frac{1}{2}\int_t^T|\alpha_n(s)|^2ds\leq a+o(1)\quad[n\to \infty].$$
Thus $(\alpha_n)_{n\in\N}$ is bounded in $L^2(t,T)$ and extracting a subsequence if necessary, one can suppose that the sequence $(\alpha_{n})_{n\in\N}$ converges weakly in $L^2(t,T)$ to some $\overline
\alpha\in L^2(t,T)$. By (A2), $(Y^{t,x,\alpha_{n}})_{n}$ converges uniformly on $[t,T]$ to
$Y^{t,x,\overline\alpha}$ and since for all $n$,
$Y^{t,x,\alpha_{n}}=g_n$ the sequence $(g_n)_{n\in\N}$ 
converges to some
$\overline g=Y^{t,x,\overline\alpha}$ in $\X$. Moreover,
\[\lambda_{t,x}(\overline g)\leq
\frac{1}{2}\int_t^T|\overline\alpha_s|^2ds \leq
\liminf_{n\to\infty}\frac{1}{2}\int_t^T|\alpha_n(s)|^2ds\leq a,\]
hence $\overline g\in\K$. We have proved that $\K$ is compact, and  the proof of point (iii) is complete.
\quad $\diamond$

\section{Proof of assertion (A1)}\label{Stube}
Fix $g_0\in \X$ and $r>0$ and consider the ball $ \B(g_0,r)$.  
The aim of this section is to prove that for each $(t,x)\in[0,T]\times \oGO$, the probability $u_\eps(t,x)$ defined by
$u_\eps(t,x)= P[X^{t,x,\eps}\in\B(g_0,r)]$  satisfies 
$\displaystyle \limsup_{\varepsilon \to 0}\left\{ -\varepsilon^2 \ln u_\eps(t,x)\right\}\le \Lambda_{t,x}(\B(g_0,r)). $

 \subsection*{Step 1. From a probability to a PDE} 
We first interpret the probability $u_\eps(t,x)$ as the value function of an optimal stopping problem.

Let us define the tube ${\bf B}$ as the set
\begin{equation}\label{B}
{\bf B}= \{(t,x)\in [0,T]\times \oGO, \; |x-g_0(t)|<r\}.
\end{equation}
\begin{proposition}\label{Ptao} $u_\eps(t,x)$ is the value function of the following optimal stopping problem
$$ u_\eps(t,x)=\inf_{\theta\in T_t} E[{\bf 1_B} (\theta, X^{t,x,\varepsilon}_\theta)],$$ where $T_t$ is the set of stopping times $\theta$ with value in $[t,T]$.
\end{proposition}
The proof can be found at the end of this section.

\vspace{1mm}

We now recall that the value function of an optimal stopping time problem is a viscosity solution of a variational inequality.

 More precisely for each bounded Borel function $\psi$ on $[0,T]\times \R^d$ consider 
\begin{equation}\label{upsi}
U_\eps[\psi](t,x)=\inf_{\theta\in T_t} E[\psi (\theta, X^{t,x,\varepsilon}_\theta)],
\end{equation}
where $ X^{t,x,\varepsilon}$ is the solution of (\ref{Eoe}), then  $U^\eps[\psi]$ is a solution of
\begin{equation}\label{IVpsi}
\left\{\begin{array}{l}
{\displaystyle  \max\left(-\dt{u}+ {\mathcal L}_\eps u,u-\psi\right)=0\mbox{ in 
}[0,T)\times\GO}\vspace{2mm}\\
{\displaystyle \dpart{u}{\gamma}=0 \mbox{ in }[0,T)\times \dGO,\
u(T)=\psi(T)\mbox{ on }\oGO.}
\end{array}\right.
\end{equation}
where $\displaystyle {\mathcal L}_\eps u=-\frac{\eps^2}{2}{\rm
Tr}[\sigma_\eps\sigma_\eps^TD^2u]-b_\eps\cdot
Du$.
\begin{proposition} Assume (\ref{gamma}), (\ref{Dom}) and (\ref{bsigma1}).
Then the function $U_\eps[\psi]$ defined by (\ref{upsi})
is a viscosity solution of (\ref{IVpsi}).
\end{proposition}
 This result is a standard consequence of the well-known Dynamic Programming Principle. Under regularity conditions the proof  goes back to \cite{BL}. For a general proof of the Dynamic Programming Principle see \cite{EKLM} or \cite{BT09}.\\ 

This gives that $u_\eps=U_\eps[{\bf 1_B}]$ is a solution of the variational inequality (\ref{IVpsi})  with obstacle ${\bf 1_B}$.

\subsection*{Step 2. The logarithmic transform} 

For all function $\psi$ nonnegative and bounded away from $0$, let $V_\eps[\psi]$ be defined by 
\begin{equation}\label{ve}V_\eps[\psi]=-\eps^2\ln(
U_\eps[\psi]).
\end{equation}
  Then $V_\eps[\psi]$ is a
viscosity solution of the following the variational inequality with obstacle $\eps^2\ln(\psi)$
\begin{equation}
\label{eqv}
\left\{\begin{array}{l}
{\displaystyle  \min\left(\!\!-\dt{V}+H_\eps(D^2V, DV),V-\eps^2\ln(\psi) \!\!\right)=0\mbox{ in 
}[0,T)\times\GO,}\vspace{2mm}\\
{\displaystyle \dpart{V}{\gamma}=0 \mbox{ in }[0,T)\times \dGO,\
V(T)=\eps^2\ln(\psi(T))\mbox{ on }\oGO,}
\end{array}\right.
\end{equation}
where $H_\eps(D^2V, DV)= -\frac{\eps^2}{2}{\rm Tr}[\sigma_\eps\sigma_\eps^TD^2V]+\frac{1}{2} |\sigma_\eps^TDV|^2\!\!-b_\eps\cdot DV$.
\\

Formaly, $v_\eps:=-\eps^2 \ln(u_\eps)=V[{\bf 1}_{{\bf B}}]$ is a viscosity solution of variational inequality (\ref{eqv}) with singular obstacle  $\chi_{\bf{B}^c}=-\eps^2 \ln( {\bf 1}_{{\bf B}})$ defined by 
$$\chi_{\bf{B}^c}(x)=  \left\{\begin{array}{ll} \infty &\mbox{ if }x\in{\bf{B}^c},\\  0&\mbox{ if }x\in{\bf{B}}.\end{array}\right.$$ 

In order to avoid the singularity, we seek  now to approximate the original obstacle ${\bf 1}_{\bf B}$ in such a way that after the logarithmic transform, the obstacle becomes $A{\bf 1}_{\bf{B}^c}$ with $A>0$. We define
for all $A,\eps>0$,  the real valued functions $\psi_\varepsilon^A$,   $u_\eps^A$ and $v_\eps^A$ by 
\begin{equation}\label{psieps}
\psi_\eps^A=\exp(-A
{\bf 1}_{{\bf{B}}^c}/\eps^2), \quad \quad  u_\eps^A=U[\psi_\eps^A]\quad \mbox{ and }\quad v_\eps^A=V[\psi_\eps^A].
\end{equation}
Note that $\psi_\eps^A\ge {\bf 1}_{\bf B}$, hence $ u_\eps^A\ge u_\eps$ and $v_\eps^A\le v_\eps$. As our aim is to majorate $\limsup v_\eps$,  it seems  at first that we have the  inequality from the wrong side. However, the following lemma shows that  we can reduce ourselves to the study of $v_\eps^A$.
\begin{lemma} \label{LuA}
For all $A>0$, and for all  $(t,x)$ $\in$  $[0,T]\times\oGO$, we have
$$\limsup_{\eps\to 0}v^A_\eps =\limsup_{\eps\to 0}v_\eps\wedge A.$$
\end{lemma}
The proof can be found at the end of this section.

\vspace{1mm}

Clearly $v_\eps^A$ is a viscosity solution of  variational inequality (\ref{eqv}) with obstacle $A{\bf 1}_{\bf{B}^c}$.
\subsection*{Step 3. Passing to the limit}

When $\varepsilon$ goes to $0$, equation (\ref{eqv}) with obstacle $A{\bf 1}_{\bf{B}^c}$ converges to the following variational inequality with obstacle $A{\bf 1}_{\bf{B}^c}$
\begin{equation}
\label{eqvl}
\left\{\begin{array}{l}
{\displaystyle  \min\left(-\dt{v}+\frac{1}{2}|\sigma^TDv|^2-b\cdot
Dv,v-A{\bf{1}}_{{\bf{B}}^c}\right)=0\mbox{ in }[0,T)\times\GO,}\vspace{2mm}\\
{\displaystyle \dpart{v}{\gamma}=0 \mbox{ in }[0,T)\times \dGO,\
v(T)=A{\bf{1}}_{{\bf{B}}^c}(T)\mbox{ on }\oGO.}
\end{array}\right.
\end{equation}

By a general stability result for viscosity solutions (see \cite{Barles94} or \cite{BP1}),
the half-relaxed upper-limit $\limsup^* v^A$  defined for all $(t,x)$ in $[0,T]\times \oGO$ by
$$\limsup\ \!\!\!^*\
v_\eps^{A}(t,x)=\limsup_{\arglim{(s,y)\to (t,x)}{ \eps\to 0}}
v_\eps^{A}(s,y). $$
is a viscosity subsolution of the limit equation (\ref{eqvl}).

\subsection*{Step 4. A first order mixed optimal control-optimal stopping problem: back to the action functional} 

 We now study a value function of a mixed optimal control-optimal stopping problem which appears to be the maximal viscosity subsolution of equation (\ref{eqvl}), and which we compare with $\Lambda_{t,x}(\B)$.

For each bounded Borel function $\psi$,  and for all $(t,x)\in [0,T]\times\oGO$, define the following value function
\begin{equation}\label{vA} v[\psi](t,x)=\inf_{\alpha\in L^2(t,T)}\sup_{\theta\in
[t,T]}\left\{\frac{1}{2}\int_t^\theta|\alpha_s|^2ds+ \psi(\theta,Y^{t,x,\alpha}_\theta)\right\}
\end{equation}
where $Y^{t,x,\alpha}$ is the unique solution of (\ref{Y}).
\begin{proposition}\label{PvA}
\begin{enumerate}
\item For each bounded Borel function $\psi$ the function $v[\psi^*]$ defined by (\ref{vA})  is  the maximal usc viscosity subsolution of the variational inequality (\ref{eqvl}) with obstacle $\psi$.
\item  For all $A>0$, $(t,x)\in [0,T]\times\oGO$, one has
$v[A{\bf 1}_{\bf{B}^c}](t,x) \le A\wedge\Lambda_{t,x}(\B)$.
\end{enumerate}
\end{proposition}
The proof can be found in Appendix C for point (1.) and at the end of this section for point (2.).

\vspace{1mm}

\subsection*{Conclusion}

By Lemma \ref{LuA}, by using the half-reaxed semi-limit method, and by Proposition \ref{PvA} (1) and (2), we have, for each $A>0$, 
$$   A\wedge\limsup_{\eps\to 0}\left\{-\eps^2\ln u_\eps(t,x)\right\}  \le \limsup{}\!^*v_\eps^A(t,x)\le v[A {\bf 1}_{\bf{B}^c}  ]\le  A\wedge\Lambda_{t,x}(\B).$$
As the inequality holds for all $A>0$, and proof of (A1) is complete.\quad $\diamond$\\

We now turn to the proofs of Propositions \ref{Ptao}, \ref{PvA} and  Lemma \ref{LuA}.\\

\noindent{\em Proof of Proposition \ref{Ptao}:}\\  Obviously,
if $X^{t,x,\eps}\in \B$ then $(\theta,X^{t,x,\eps}_\theta)\in
{\bf B}$ for all $\theta\in T_t$,
hence
$\ue(t,x)\leq 
E[{\bf 1_B}(\theta,X^{t,x,\eps}_\theta)]$. Taking the infimum over $\theta \in T_t$ we obtain
$$u_\eps(t,x)\leq \inf_{\theta\in T_t} E[{\bf 1_B} (\theta, X^{t,x,\varepsilon}_\theta)].$$
Conversely, let $\tilde\theta$ be the first exit time of $(s,X^{t,x,\eps}_s)_{s\in[t,T]}$ from  ${\bf B}$ 
$$\tilde\theta=\inf\{s\geq t, \;|X^{t,x,\eps}_s-g_0(s)|\geq r\}.$$
Then $\tilde \theta\wedge T$ is a stopping time which takes value in $[t,T]$. Let us show that almost surely 
$$\{ (\tilde \theta\wedge T, X^{t,x,\eps}_{\tilde\theta\wedge T})\in
{\bf B}\}\subset  \{X^{t,x,\eps}\in \B\}. $$
Indeed, suppose   $ (\tilde \theta(\omega)\wedge T, X^{t,x,\eps}_{\tilde\theta(\omega)\wedge T})\in
{\bf B}$, then  $\tilde\theta(\omega)>T$. Hence for all $s\in[t,T]$ we have
$|X^{t,x,\eps}_s(\omega)-g_0(s)|<r$,  which means, as both $X^{t,x,\eps}_.(\omega)$ and $g_0(.)$
 are continuous on $[t,T]$ that $\|X^{t,x,\eps}(\omega)-g_0\|_\infty<r$, hence $X^{t,x,\eps}(\omega)\in \B$ and the proof is complete.\quad $\diamond$\\

\noindent{\em Proof of Lemma \ref{LuA}: }\\
 Fix $(t,x)$ $\in$  $[0,T]\times \oGO$, and $A>0$. Clearly
$e^{-A/\eps^2}\vee {\bf 1_{B}}(t,x)\le \psi_\eps^A(t,x)\leq
{\bf 1_{B}}(t,x)+e^{-A/\eps^2}.$
This gives easily,
$$
e^{-A/\eps^2}\vee u_\eps(t,x)\le u_\eps^A(t,x)\le  u_\eps(t,x)+e^{-A/\eps^2}.
$$
As for any nonegative sequence  $ (u_\eps)$ one has\\  
$\displaystyle \limsup_{\eps\to
0}\left\{-\eps^2\ln(u_\eps+e^{-A/\eps^2})\right\}=A\wedge\limsup_{\eps\to 
0}\left\{-\eps^2\ln u_\eps\right\},$
we obtain 
$$ \displaystyle A\wedge \limsup_{\eps\to0}\left\{-\eps^2\ln u_\eps\right\}\ge \limsup_{\eps\to 
0}\left\{-\eps^2\ln u^A_\eps\right\}\ge A\wedge\limsup_{\eps\to 
0}\left\{-\eps^2\ln u_\eps\right\},
$$
which completes the proof of the lemma.\quad $\diamond$\\

\noindent{\em Proof of Proposition \ref{PvA}:} \\
Point (1) is detailed in Appendix C  (Proposition \ref{Puu}).

Let us prove now the second point.
Obviously, $v^{A}(t,x)\leq A$.
Now, if $\Lambda_{t,x}(\B)<A$,  for each $\eta>0$ such that $\Lambda_{t,x}(\B)+\eta <A$ there exists $\tilde \alpha$ $\in$ $L^2(t,T)$ such that $Y^{t,x,\tilde \alpha}$ $\in $ $\B$ and 
$$ \Lambda_{t,x}(\B)\leq \frac{1}{2} \int_t^T
|\tilde \alpha_s|^2ds \leq \Lambda_{t,x}(\B)+\eta.$$ Thus, for any $\theta\in
[t,T]$, one has ${\bf{1}}_{{\bf{B}}^c}(\theta,Y^{t,x,\alpha}_\theta)=0$ so that
\[v^A(t,x)= \inf_{\alpha\in L^2(t,T)}\sup_{\theta\in[t,T]}\frac{1}{2}\int_t^T
|\alpha_s|^2ds \leq \frac{1}{2}\int_t^\theta
|\tilde \alpha_s|^2ds\leq  \Lambda_{t,x}(\B)+\eta.\]
We have proved that $ v^A(t,x) \le A\wedge\Lambda_{t,x}(\B)$.\quad   $\diamond$

\section{Proof of assertion (A2)}\label{Stubes}
Let $(g_n)_{n\in \N}$ be a sequence of functions in $\X$ and $(r_n)_{n\in \N}$, 
a sequence of positive reals. For each nonempty finit subset  $I$ of $\N$ and for all $(t,x)\in [0,T]\times \oGO$, we define 
\begin{equation}
\label{tueps}
u^I_\eps(t,x)= P[X^{t,x,\eps}\in\G], \quad \mbox{ where } \quad \G =\bigcap_{i\in I}\B(g_i,r_i)^c,
\end{equation}
and we prove that $\displaystyle \liminf_{\eps\to 0} \left\{-\eps^2 \ln u^I_\eps(t,x) \right\}\ge \Lambda_{t,x}(\G)$.\\

 In the following we will denote by $\theta_I$ a multiple stopping time $(\theta_i)_{i\in I}$ with $\theta_i\in T_t$ for each $i\in I$, and we write $\theta_I\in T_t^I$.
\subsection*{Step 1. From a probability to a PDE}
 We first interpret $u^{I}_\eps(t,x)$ as the value of an optimal multiple stopping times problem.
\begin{proposition}
For all $(t,x)\in [0,T]\times\oGO$, for all $\eps>0$, and for each finit nonempty subset $I$ of $\N$ 
\[u^I_\eps(t,x)=\sup_{\theta_I\in T^I_t}E\left[\prod_{i \in I}({\bf 1}_{{\bf B}^c_i}(\theta_i,X^{t,x,\eps}_{\theta_i}))\right]\]
where for all $i\in I$, ${\bf B}_i=\{(t,x)\in[0,T]\times \oGO ; \; |x-g_i(t)|<r_i\}$.
\end{proposition}

\noindent{\em Proof: } It is similar to the proof of Proposition \ref{Ptao}. The 
main
difference is in the choice of the optimal stopping time which is here
 $\tilde\theta_I\in T_t^I$ where for $i\in I$, $\tilde\theta_i$ is the first exit time in $[t,T]$ of 
 $(s,X^{t,x,\eps}_s)$ from ${\bf B}_i$.\quad $\diamond$\\

  As by the logarithmic transform the obstacles would take there values in $[0,\infty]$, 
 the first task is to approximate ${\bf 1}_{{\bf B}^c_i}$.
 For all $A,\eps>0$ and for all $i\in\N$, we define
\begin{equation}\label{psiIAe} 
\psi_\eps^{\{i\},A}=\exp\left(-\frac{A{\bf{1}}_{{\bf{B}}_i}}{\eps^2}\right).
 \end{equation} Note that 
for all $\eps>0$  we have 
 ${\bf 1}_{{\bf B}^c_i}\leq \psi_\eps^{i,A}.$\\
 We define for each nonempty finit subset $I$ of $\N$, for each $A,\eps>0$, and each $(t,x)\in[0,T]\times \oGO$ the value function 
\begin{equation}\label{uIAe} 
  u^{I,A}_\eps(t,x)=\sup_{\theta_I\in T_t^I}
 E\left[\prod_{i\in I}\psi_\eps^{\{i\},A}(\theta_i,X^{t,x,\eps}_{\theta_i})
  \right].
 \end{equation}
 Clearly
 \begin{equation}
 \label{uI-uIA}
 u^I_\eps(t,x)\leq u^{I,A}_\eps(t,x).
 \end{equation}
 
 We now proceed to the reduction of the multiple stopping problem to a single stopping problem. 
More precisely let us define  for all finit subset $I$ of $\N$ containing two or more elements, for all $A, \eps>0$, and for all  $(t,x)$ $\in$ $[0,T]\times \oGO$,  
\[\psi_\eps^{I,A}(t,x)=\max_{i\in I}\left\{\psi_\eps^{\{i\},A}(t,x)
 u_\eps^{I\backslash \{i\},A}(t,x)\right\}.\]
One has by Theorem 3.1 in \cite{TAM10} 
 \[  u_\eps^{I,A}(t,x) =  \sup_{\theta\in T_t}
 E\left[\psi_\eps^{I,A}(\theta,X^{t,x,\eps}_\theta)\right].\]
Now, one can show (cf \cite{BT09}) that  $ u^{I,A}_\eps$ is a viscosity subsolution of the following
variational inequality
\[\left\{\begin{array}{l}
{\displaystyle  \min\left(-\dt{u}+{\mathcal L}_\eps u ,u-\psi_\eps^{I,A}\right)=0\mbox{ in 
}[0,T)\times\GO},\vspace{2mm}\\
{\displaystyle \dpart{u}{\gamma}=0 \mbox{ in }[0,T)\times 
\dGO,\
u(T)=\psi_\eps^{I,A}(T)
\mbox{ on }\oGO.}
\end{array}\right.\]

\subsection*{Step 2. The logarithmic transform}

For all nonempty finit subset $I$ of $\N$, for all $ A,\eps>0$, let 
$v_\eps^{I,A}$ be defined on $ [0,T]\times \oGO$ by
\begin{equation}\label{vIedA}
v_\eps^{I,A}=-\eps^2\ln  u_\eps^{I,A}.
\end{equation}
Then $v_\eps^{I,A}$ is a viscosity supersolution  of the following variational inequality
\[\left\{\begin{array}{l}
{\displaystyle  \max\left(-\dt{v}+H_\eps(D^2v,Dv), v- \phi_\eps^{I,A}\right)=0\mbox{ in 
}[0,T)\times\GO}\vspace{2mm}\\
{\displaystyle \dpart{ v}{\gamma}=0 \mbox{ in }[0,T)\times 
\dGO,\
 v(T)= \phi_\eps^{I,A}(T)\mbox{ on }\oGO}
\end{array}\right.\]
where, for all nonempty finit subset $I$ of $\N$, for all $A, \eps>0$,
\begin{equation}\label{phiIA2}
\phi_\eps^{I,A}=\left\{\begin{array}{ll} A{\bf{1}}_{{\bf{B}}_i}& \mbox{ if } I=\{i\} \mbox{ with } i\in \N,\\
\displaystyle \min_{i\in I}\left\{A{\bf{1}}_{{\bf{B}}_i}+
v_\eps^{I\backslash\{i\},A}\right\}& \mbox{ if } {\rm card }\, I \ge 2.
\end{array}\right.
\end{equation}

\subsection*{Step 3. A mixed optimal control-optimal multiple stopping problem} Let us turn now to the study of a mixed optimal control-optimal multiple stopping problem. The value function of this problem will be shown  to be smaller than the  half-relaxed lower limit  $\lim_* v_\eps^{I,A}(t,x)$ and greater than $\Lambda_{t,x}(\G)\wedge A$.

For all finit and nonempty subset $I$ of $\N$  and for all $(t,x)\in [0,T]\times\oGO$, define the following value function
\begin{equation}\label{vIA} v^{I,A}(t,x)=\inf_{\alpha\in L^2}\inf_{\theta_I\in
[t,T]^I}\left\{\frac{1}{2}\int_t^{\vee_{i\in I}\theta_i}|\alpha_s|^2ds
+\sum_{i\in I}A{\bf{1}}_{{\bf{B}}_i}(\theta_i,Y^{t,x,\alpha}_{\theta_i})\right\},
\end{equation}
where $Y^{t,x,\alpha}$ is the unique solution of (\ref{Y}).\\

This mixed  optimal {\em multiple} stopping problem can be reduced to a mixed optimal {\em single} stopping problem. More precisely, consider for each bounded real valued measurable $\phi$ defined on $[0,T]\times \oGO$ and for each $(t,x)\in[0,T]\times \oGO$ the following value function 
 \begin{equation}\label{vphi}
 v[\phi](t,x)=\inf_{\alpha\in L^2(t,T)}\inf_{\theta\in[t,T]}\left\{
\frac{1}{2}\int_t^\theta
|\alpha_s|^2ds+\phi(\theta,Y^{t,x,\alpha}_\theta)\right\},
\end{equation}
where $Y^{t,x,\alpha}$ is the unique solution of (\ref{Y}).\\
Define also for all nonempty finit subset $I$ of $\N$, for all $A>0$,
\begin{equation}\label{phiIA}
\phi^{I,A}=\left\{\begin{array}{ll} A{\bf{1}}_{{\bf{B}}_i}& \mbox{ if } I=\{i\} \mbox{ with } i\in \N,\\
\displaystyle \min_{i\in I}\left\{A{\bf{1}}_{{\bf{B}}_i}+
v^{I\backslash\{i\},A}\right\}& \mbox{ if } {\rm card\; } I \ge 2.
\end{array}\right.
\end{equation}
\begin{proposition}\label{PvIA}
Let $I$ be a finit subset of $\N$ and $A>0$, and consider the function
$v^{I,A}$ defined by (\ref{vIA}).
Then 
\begin{enumerate}
\item for each $(t,x)\in [0,T]\times \oGO$ one has $v^{I,A}(t,x)= v[\phi^{I,A}](t,x)$ where $\phi^{I,A}$ is defined by
(\ref{phiIA}),
\item  one has for all  $(t,x)\in [0,T]\times\oGO$,
$\displaystyle v^{I,A}(t,x) \ge A\wedge\Lambda_{t,x}(\G).$
\end{enumerate}
\end{proposition}

\noindent{\em Proof: } The proof of (1) is the concequence of a reduction result for optimal multiple stopping problems. It is detailed in Appendix C (Proposition \ref{Preduction}).

\noindent Let us prove (2). Suppose $v^{I,A}(t,x)<A$, then for each $\eta>0$ such that $v^{I,A}(t,x)+\eta <A$, there exists
$\theta_I\in [0,T]^N$ and $\alpha \in L^2(t,T)$  such that
\[ \frac{1}{2}\int_t^{\vee_{i\in I}\theta_i}|\alpha_s|^2ds +\sum_{i\in I}
A{\bf{1}}_{{\bf{B}}_i}(\theta_i,Y^{t,x,\alpha}_{\theta_i})\le v^{I,A}(t,x)+\eta < A.\] This means in particular that
$\displaystyle   (\theta_i,Y^{t,x,\alpha}_{\theta_i})\in  {\bf{B}}_i^c $
for all $ i\in I$ and therefore $Y^{t,x,\alpha}\in
\G$. We set, for all $s\in [t,T]$,
\[\tilde \alpha_s=\left\{\begin{array}{l}
\alpha_s\mbox{ if }s\leq \vee_{i\in I}\theta_i\\
0\mbox{ otherwise}.
\end{array}\right.\]
Then again $Y^{t,x,\tilde \alpha}\in\G$ and
\[\Lambda_{t,x}(\G)\leq \frac{1}{2}\int_t^T|\tilde \alpha_s|^2ds \le  v^{I,A}(t,x)+\eta\]
and letting $\eta$ to 0 the proof is complete.\quad $\diamond$

\vspace{2mm}

We now give some results concerning the mixed optimal single stopping problem (\ref{vphi}), and its links with the following  variational inequality:
\begin{equation}\label{eqvphi}\left\{\begin{array}{l}
{\displaystyle  \max\left(-\dt{v}+\frac{1}{2}|\sigma^TD v|^2-b\cdot
D v, v- \phi\right)=0\mbox{ in 
}[0,T)\times\GO}\vspace{2mm}\\
{\displaystyle \dpart{ v}{\gamma}=0 \mbox{ in }[0,T)\times 
\dGO,\
 v(T)= \phi(T)\mbox{ on }\oGO}
\end{array}\right.
\end{equation}

\begin{proposition}\label{Pvphi}
 Let $\phi:[0,T]\times \oGO$ be a measurable, bounded, real valued function then the function $v[\phi_*]$ is the minimal  lsc viscosity supersolution of the  variational inequality (\ref{eqvphi}).
\end{proposition}
 The proof is similar  and even simpler than the proof of Proposition \ref{Puu}. Let us remark that
this result is well known for deterministic systems with Lipschitz coefficients in $\R^n$ (see Barles and Perthame \cite{BP1}). The main difficulty in the present case is to prove the minimality of $v[\phi_*]$. This point is the concequence of a strong comparison result for equation (\ref{eqvphi}) when the obstacle is bounded and continuous on $[0,T]\times \oGO$. The proof of this strong comparison result, which is highly technical, is detailed in Appendix D.

\subsection*{Step 4. Passing to the limit}

Let us now prove the following lemma
\begin{lemma}\label{LvIA} For each finit non-empty subset $I$ of $\N$ and for each $A>0$, one has
$$\liminf\,\!\!_* \,v^{I,A}_\eps(t,x)\ge v[\phi^{I,A}](t,x).$$
\end{lemma}
\noindent{\em Proof: } The result is established by induction on the cardinal of $I$. If $I=\{i\}$ for some $i$ in $\N$, then $\liminf_* \,v^{\{i\},A}_\eps$ is a viscosity supersolution of (\ref{eqvphi}) with $\phi=\phi^{\{i\},A}$. By Proposition \ref{Pvphi}, $v[\phi^{\{i\},A}]$ is the minimal viscosity supersolution of the same equation, hence the proof is complete.

Suppose now that $I$ has $N$ elements with $N\geq 2$ and that the lemma holds for any subset $J$ of $\N^*$ containing $N-1$ elements. 
Then by using the induction hypothesis on formula (\ref{phiIA}), one  has $\liminf_* \phi^{I,A}_\eps= \tilde \phi^{I,A} \ge  \phi^ {I,A}$. By a stability result 
$\liminf_*v^{I,A}_\eps$  is a viscosity supersolution of  (\ref{eqvphi}) with obstacle $ \tilde \phi^{I,A}$. By Proposition \ref{Pvphi}, as $\tilde \phi^{I,A}$ is lsc, the minimal viscosity supersolution of this equation is $ v[\tilde \phi^{I,A}]$. Now as $\tilde \phi^{I,A} \ge  \phi^ {I,A}$ one clearly  has $v[\tilde 
\phi^{I,A}] \geq v[\phi^{I,A}]$.

Finally we have $\liminf_*v^{I,A}_\eps\ge v[\tilde \phi^{I,A}]\ge v^{I,A} $ which completes the proof of the lemma. \quad $\diamond$

\subsection*{Conclusion}
For  all finit nonempty subset $I$ of $\N$, for all $A>0$, and for all $(t,x)$ $\in$  $[0,T]\times\oGO$ one has, by inequality (\ref{uI-uIA}),
$ u^I_\eps(t,x)\leq  u_\eps^{I,A}(t,x)$, hence, by Lemma \ref{LvIA} and by Proposition \ref{PvIA},
\[ \liminf_{\eps\to 0}\left\{-\eps^2\ln u^I_\eps(t,x)\right\}\ge \liminf\,\!\!_* \,v^{I,A}_\eps(t,x)\ge v[\phi^{I,A}](t,x) \ge \Lambda_{t,x}(\G)\wedge A. \]
 The proof of (A3) is complete. \quad $\diamond$

\section{Appendix}

\section*{ Appendix A: the test-function}
\begin{lemma}\label{LT} We assume that  $\gamma$ and $\GO$ satisfy
(\ref{gamma}) and (\ref{Dom}). Then,  for all
$\eps,\rho >0$, there exists $\psi_{\eps,\rho}\in 
C^1(\oGO\times\oGO,\R)$ such that,
\[\begin{array}{ll}
(\psi i) & \dis \forall x,y\in \oGO,\quad\myfrac{1}{2} 
\myfrac{|x-y|^2}{\eps^2} -K
\myfrac{\rho^2}{\eps^2}\leq
\psi_{\eps,\rho}(x,y) \leq K ( \myfrac{|x-y|^2}{\eps^2} 
+\myfrac{\rho^2}{\eps^2}), \vspace{2mm}\\
(\psi ii) & \left\{\begin{array}{l}
\dis \forall x,y\in \oGO,\quad|D_x\psi_{\eps,\rho}(x,y)+D_y 
\psi_{\eps,\rho}(x,y)|
\leq K ( \myfrac{|x-y|^2}{\eps^2} +\myfrac{\rho^2}{\eps^2}) 
,\vspace{2mm}\\
\dis \forall x,y\in \oGO,\quad |D_x\psi_{\eps,\rho}(x,y)|+|D_y 
\psi_{\eps,\rho}(x,y)|
\leq K ( \myfrac{|x-y|}{\eps^2} +\myfrac{\rho^2}{\eps^2}),
\end{array}\right.\vspace{2mm}\\
(\psi iii) & \left\{\begin{array}{l}
\dis \forall x \in \dGO, y\in \oGO, \quad 
D_x\psi_{\eps,\rho}(x,y)\cdot \gamma(x)>0, \vspace{2mm}\\
\dis \forall y \in \dGO,x\in \oGO, \quad 
D_y\psi_{\eps,\rho}(x,y)\cdot \gamma(y)>0.
\end{array}\right.
\end{array}\]
for some constant  $K$ depending only on $\GO$,
$||\gamma||_{\infty}$ $||\gamma||_{Lip}$ and $c_0$.
\end{lemma}

We use ideas from \cite{Barles93}.

\noindent{\em Proof: } We first define the Lipschitz continuous $\R^d$-valued function
$\mu$ on $\dGO$ by
\[\mu(x)=\frac{\gamma(x)}{\gamma(x)\cdot n(x)}\]
as well as its smooth approximation $(\mu_{\rho})_{\rho>0}$ such that 
for all $\rho>0$,
$\mu_{\rho}\in C^1(\R^d,\R^d)$,
\[\|\mu_\rho\|+\|D\mu_\rho\|\leq K_1\]
for some constant $K_1>0$ independent of $\rho$ and for all $x\in\dGO$
\[|\mu_{\rho}(x)-\mu(x)|\leq \rho.\]

Then we set,
\begin{eqnarray*}
\phi_{\eps,\rho}(x,y)&=&    \myfrac{|x-y|^2}{\eps^2} + 2 \frac{(x-y) 
}{\eps}\cdot
 \mu_\rho(\frac{x+y}{2})  \frac{(d(x)-d(y))}{\eps}  +A \frac{ 
(d(x)-d(y))^2}{\eps^2}.
\end{eqnarray*}
We can choose the
constant $A>0$ large enough in order to get, for some constant $K_2>0$
and for all $\eps$, $\rho>0$,
$$\begin{array}{ll}
(\phi i) & \dis \forall x,y\in \oGO,
\quad\frac{1}{2} \myfrac{|x-y|^2}{\eps^2}\leq \phi_{\eps,\rho}(x,y) 
\leq K_2  \myfrac{|x-y|^2}{\eps^2} ,\vspace{2mm}\\
(\phi ii) & \left\{\begin{array}{l}
\dis \forall x,y\in \oGO,\quad|D_x\phi_{\eps,\rho}(x,y)+D_y 
\phi_{\eps,\rho}(x,y)|\leq K_2  \myfrac{|x-y|^2}{\eps^2}
,\vspace{2mm}\\
\dis \forall x,y\in \oGO,\quad |D_x\phi_{\eps,\rho}(x,y)|+|D_y 
\phi_{\eps,\rho}(x,y)| \leq K_2  \myfrac{|x-y|}{\eps^2} ,
\end{array}\right.\vspace{2mm}\\
(\phi iii) & \left\{\begin{array}{l}
\dis \forall x \in \dGO, y\in \oGO, \quad 
D_x\phi_{\eps,\rho}(x,y)\cdot
\gamma(x)\geq -K_2[\frac{|x-y|^2}{\eps^2}+ \frac{\rho^2}{\eps^2}],
\vspace{2mm}\\
\dis \forall y \in \dGO,x\in \oGO, \quad D_y\phi_{\eps,\rho}(x,y)\cdot
\gamma(y)\geq -K_2[\frac{|x-y|^2}{\eps^2}+ \frac{\rho^2}{\eps^2}].
\end{array}\right.
\end{array}$$

Indeed ($\phi i$) comes from a simple application of Cauchy-Schwarz inequality, and from the fact that $d$ is Lipschitz continuous. Now for all $U=(u,v)\in \R^d\times \R^d$, we have

\begin{eqnarray*}D\phi_{\eps,\rho}(x,y).U &= &\frac{2(x-y)\cdot(u-v)}{\eps^2}+ 2\frac{u-v}{\eps^ 2}\cdot \mu_{\rho}(\frac{x+y}{2})(d(x)-d(y))\\
&&- 2\frac{x-y}{\eps^ 2}\cdot\mu_{\rho}(\frac{x+y}{2})(n(x)u-n(y)v)\\
&&+\frac{(x-y)}{\eps^ 2}\cdot D\mu_{\rho}(\frac{x+y}{2})(u+v)(d(x)-d(y))\\
&&- \frac{2A}{\eps^2}(d(x)-d(y))(n(x)u-n(y)v).
\end{eqnarray*}
Taking $U=(u,u)$, as both $d$ and $n$ are Lipschitz continuous, and using Cauchy-Schwarz inequality we obtain straightforwardly the first inequality in ($\phi ii$). The second inequality in  ($\phi ii$) is clear. 

Let us  now  prove ($\phi iii$). By symmetry, there is only one inequality to prove. Take $x\in \dGO$ and $y\in \oGO$, and take $U=(\gamma(x),0)$, and recall that $\gamma(x).n(x)\ge c_0>0$.
The sum of all  the terms that have $(d(x)-d(y))=-d(y)$ can be made nonnegative for $A$ large enough.
The remaining term is, taking $\displaystyle \frac{2(x-y)}{\eps^2}$ in factor,
$$ \gamma(x)- \mu_{\rho}(\frac{x+y}{2}) n(x)\cdot\gamma(x)= \left(\mu(x)-\mu_{\rho}(\frac{x+y}{2})\right)(n(x)\cdot\gamma(x)).$$
Writing \\

$|\mu(x)-\mu_{\rho}(\frac{x+y}{2})|\le | (\mu(x)-\mu(\frac{x+y}{2})|+ |\mu(\frac{x+y}{2})-\mu_{\rho}(\frac{x+y}{2})|\le  K |\frac{x-y}{2}|+\rho,$ \\

\noindent we have completed the proof of ($\phi iii$).

\vline 

Finally, we set, for $x,y\in \oGO$,
\begin{equation}
   \psi_{\eps,\rho}(x,y) = e^{C(2\|d\|_\infty-d(x)-d(y))}
      \phi_{\eps,\rho}(x,y) -B\frac{\rho^2}{\eps^2}(d(x)+d(y)).
\end{equation}
By choosing $B$, then $C$ large enough, we obtain the desired result. \quad $\diamond$

\section*{ Appendix B. A deterministic reflection problem}

In this section, we suppose that $b$ and $\sigma$ satisfy (\ref{bsigma1}). 
We  consider for  each fixed  $\alpha\in L^2(0,T; \R^m)$ and for each fixed $(t,x)\in[0,T]\times \oGO$  the deterministic equation with oblique reflection
\begin{equation}
\label{eqdet}
\left\{\begin{array}{l}
dY_s=\left(b(s,Y_s)-\sigma(s,Y_s)\alpha_s\right) ds-dz_s,  \quad Y_s\in\oGO, \quad \mbox{ for all
} s\in [0,T]\vspace{2mm}\\
dz_s={\bf 
1}_{\dGO}(Y_s)\gamma(Y_s)d|z|_s \quad \mbox{ for all
}s\in [0,T], \ \vspace{2mm}\\
Y_t=x \mbox{ and } z_t=0.
\end{array}\right.
\end{equation}
A solution of equation (\ref{eqdet}) is a couple $(Y,z)$ of continuous functions defined on $[0,T]$ with values in $\R^d$ such that $z$ has bounded variations, and $|z|_s$ denotes the total variation of $z$ on the interval $[0,s]$.

\begin{theorem}\label{leY}
Assume (\ref{gamma})-(\ref{Dom})-(\ref{bsigma1}) and let $\alpha\in L^2(0,T; \R^d)$.
Then, 
\begin{enumerate}
\item there exists a unique solution $(Y^{t,x,\alpha},z^{t,x,\alpha})$ of (\ref{eqdet}),
\item for each $s\in[0,T]$,  the function  $(t,x)\mapsto Y^{t,x,\alpha}_s$ is continuous and 
$$|Y_s^{t,x,\alpha}-Y_s^{t',x',\alpha}|^2\leq K(|x-x'|^2+|t-t'|^{1/4}),$$
\item for each $(t,x)\in[0,T]\times \oGO$ the function $s\mapsto Y^{t,x,\alpha}_s$ is H\"older continuous and for each $s,s'\in[0,T]$,
$$|Y^{t,x,\alpha}_s-Y^{t,x,\alpha}_{s'}|\leq K|s-s'|^{1/2}, $$
\item Assumption (A3) holds true, that is: for all $\alpha^n,\alpha\in L^2(0,T; \R^d)$, if $\alpha^n\rightharpoonup\alpha$ weakly in $L^2$, then $\|Y^{t,x,\alpha_n}-(Y^{t,x,\alpha}\|_X\to 0$
\end{enumerate}
where the constant $K$ in $\mathrm{(2)}$ and $\mathrm{(3)}$ depends  only 
on  $\GO$, $\gamma$, $c_0$ the Lipschitz constant $K_T$ of $b$, $\sigma$ and $\|\alpha\|_{L^2}$.
\end{theorem}

\subsection*{ Proof of (1)} For the sake of completness, and as the hypothesis on the coefficient $c= b-\sigma\alpha$ are slightly more general than in \cite{LS84} or \cite{DI93}, we present a complete proof. 
To that end, we use the Skorokhod problem. 
More precisely, fix  $(t,x)\in[0,T]\times \oGO$, and $\alpha\in  L^2(0,T;\R^d)$ and define   $c_t(x)= b(t,x)-\sigma(t,x)\alpha_t$ for all $(t,x)\in[0,T]\times \oGO$.

By \cite{LS84},
for  each $\X$ $\in$ $ \C([t,T];\oGO)$, there exists at least one solution 
 $(Y,z)$  of the following Skorokhod problem:
\begin{equation}\label{Sdet}
\left\{ \begin{array}{l}
\displaystyle{ Y_s=x+\int_t^s c_u(X_u)du -z_s, \quad Y\in\C([t,T];\oGO)},\vspace{2mm}\\
\displaystyle{z_s=\int_t^s{\bf 1} _{\dGO}(Y_u)\gamma(Y_u)d|z|_u, \quad z\in \C_{bv}([t,T];\oGO)}.
\end{array}\right.
\end{equation}
We next show that the solution of (\ref{Sdet}) is unique and then we prove the existence and uniqueness for the solutions of equation (\ref{eqdet}) by a fixed point argument. 
Note that, in view of the first equation of (\ref{Sdet}) it is enough to prove the uniqueness for $Y$ only. 

\vline

We define the binary relation $\mathcal{ S}$ on  $ \C([t,T];\oGO)$ in the following way: for all $ X,Y\in \C([t,T];\oGO )$, 
$Y \, \mathcal{S} \,  X$ if  and only if   
    there exists $ z\in \C_{bv}([t,T];\oGO)$ such that $(Y,z)$ is the solution of  (\ref{Sdet}).

\begin{lemma}\label{Lc} Let  $X,X',Y,Y'$ $\in$ $ \C([t,T];\oGO)$ and suppose
$Y \,\mathcal{S}\, X$ and $Y' \,\mathcal{S}\, X'$. 

Then, there
 exists $\eta>0$, depending
only on $\oGO$, $\gamma$, $c_0$ and $a_{T}$ (and independent of $t$), such that for all $u\in[t, t+\eta]\cap[0,T]$,
\[|Y_u- Y'_u|\leq \frac{1}{2}|X_u- X'_u|.\]
\end{lemma}

 \noindent{\em Proof: } We use the function $\psi_{\eps,\rho}$ defined in Lemma
 \ref{LT} with $\eps=1$, and fix $s\in[t,T]$ and we put $f_x=D_x\psi_{1,\rho}$, $f_y=D_y\psi_{1,\rho}$. We have
$$\psi_{1,\rho}(Y_s, Y'_s) =\psi_{1,\rho}(x,x)+\int_t^s f_x(Y_u, Y'_u)c_u(X_u)du+\int_t^s f_y(Y_u, Y'_u)c_u(X'_u)du$$
$$ - \int_t^s f_x(Y_u, Y'_u){\bf 1} _{\dGO}(Y_u)\gamma(Y_u)d|z|_u-\int_t^sf_y(Y_u, Y'_u){\bf 1} _{\dGO}(Y'_u)\gamma(Y'_u)d|z'|_u.$$
By $(\psi iii)$, the two last integrals of the right hand side of the above inequality are non positive. 
Write the first term of the right hand side of the previous inequality as
 $$\int_t^s( f_x+f_y)(Y_u, Y'_u)c_u(X_u)du+ \int_t^s f_y(Y_u, Y'_u)(c_u(X'_u)-c_u(X_u))du.  $$
Put  $a_u=1+| \alpha_u|$, by
 using $(\psi i)$, $(\psi ii)$ and as $|c_u(X_u)|\le K a_u$ and $|c_u(X_u)-c_u(X'_u)|\le K a_u|X_u-X'_u|$, 
 $$\hspace{-1cm}\frac{1}{2}|Y_s-Y'_s|^2\leq  K \left ( 
2\rho^2+\int_{t}^s(|Y_u- Y'_u|^2+\rho^2)a_udu\right.$$
$$\hspace{3cm}\left.+\int_t^s
 (|Y_u- Y'_u|+\rho^2)|X_u- X'_u|a_udu\right ).$$
 This equality holds independently of $\rho>0$ and its right-hand term
 is nondecreasing with $s$ therefore, by letting $\rho$ to $0$ we have for all $s\in [t,T]$, and writing for $g\in \X$, 
 $|g|_{[t,s]}=\sup\{|g(u)|, u\in[t,s]\}$, 
 $$|Y- Y'|_{[t,s]}^2\leq  K\left(\int_{t}^s a_udu\right) \left (
|Y- Y'|_{[t,s]}^2+|Y- Y'|_{[t,s]}|X- X'|_{[t,s]}\right ).$$
 Now, by using Cauchy-Shwarz inequality,
 we have for $|t-s|$ small enough,
 $$|Y-Y'|_{[t,s]}\leq \frac{2 K
 |s-t|^\frac{1}{2}\|a\|_{L^2}}{1-2 K|s-t|^{\frac{1}{2}}
 \|a\|_{L^2}}|X-X'|_{[t,s]},$$
and we can chose $\eta>0$ independently of $t$ such that 
  \[\sup_{u\in[t,t+\eta]\cap[0,T]}|Y_u-Y'_u|\le \frac{1}{2} \sup_{u\in[t,t+\eta]\cap[0,T]}|X_u-X'_u|,\] 
 which completes the proof  of the lemma.  \quad $\diamond$

\vline

Lemma \ref{Lc} shows first that for each $X$ there exists a unique $Y$ such that $X\, \mathcal{S}\, Y$. Changing notation we have proved that $\mathcal{S}: X\mapsto Y$ is a map. Lemma \ref{Lc} shows also that there exists $\eta$ such that for all $t\in[0,T]$,  $\mathcal{S}$ contracts $ \C([t,t+\eta]\cap [0,T];\oGO)$ onto itself. This  gives existence and uniqueness for (\ref{eqdet}) by a fixed point argument hence it proves  assertion (1) of the theorem.

\subsection*{Proof of  (2)} 
Let us first establish the following lemma.

\begin{lemma}\label{Li}
Fix  $\alpha,\alpha'\in L^2(0,T;\R^d)$ and $t,t'\in[0,T]$ with $t'\le t$, $x$,
$x'\in \oGO$ and define $Y=Y^{t,x,\alpha}$ and $Y'=Y^{t,x',\alpha'}$. Then, there exists a constant $K>0$ (which only depends on 
$\|\alpha\|_{L^2}$, $c_0$, $\|\gamma\|_{Lip}$, $K_T$ and $\GO$)
such that, for all $\rho>0$ and for all $s\in [t,T]$, we have
\[|Y_s-Y'_s|^2\leq K\left( 
g_\rho(s)+\int_t^s g_\rho(u) \cdot (1+|\alpha_u|+|\alpha'_u|)du\right) \]
where
\[g_\rho(s)=\rho^2+|x-x'|^2 + |t-t'|^{1/2}+ \left|\int_t^sD_{x}\psi_{1,\rho}(Y_u,Y'_u)\sigma(u, Y_u)(\alpha_u- \alpha'_u)du\right|.\]
\end{lemma}

\noindent{\em Proof: }  We put  for $(s,y)\in[0,T]\times \oGO$, $c_s(y)= b(s,y)-\sigma(s,y) \alpha_s$, $c'_s(y)= b(s,y)-\sigma(s,y) \alpha'_s$, and  $a_s= 1+|\alpha_s|$, $a'_s=1+|\alpha'_s|$. Note that there exists a constant $K$ such that for all $(s,y)\in[0,T]\times \oGO$, one has $|c_s(y)|\le K a_s$ and $|c'_s(y)|\le K a'_s$. Define as before $f_x=D_x\psi_{1,\rho}$ and $f_y=D_y\psi_{1,\rho}$. Recall also that the function $\psi_{1,\rho}$ is a continuous function from $\R^d\times \R^d$ to 
$\R^d$ which is bounded on $\oGO$ $\times$ $\oGO$ independently of $\rho\in (0,1)$, and  write
\begin{eqnarray*}
\psi_{1,\rho}(Y_s, Y'_s) &=&\psi_{1,\rho}(x,x')+\int_{t'}^t  f_y(Y_u, Y'_u)c'_u(Y'_u)du \\&& +\int_{t}^s (f_x+ f_y)(Y_u, Y'_u)c_u(X_u)du \\
&& + \int_{t}^s f_y(Y_u, Y'_u)(c'_u(Y'_u)-c_u(Y_u))du  \\
 &&- \int_{t'}^s f_x(Y_u, Y'_u){\bf 1} _{\dGO}(Y_u)\gamma(Y_u)d|z|_u\\&&-\int_{t'}^sf_y(Y_u, Y'_u){\bf 1} _{\dGO}(Y'_u)\gamma(Y'_u)d|z'|_u
\end{eqnarray*}
We then follow similar calculations as Lemma \ref{Lc}. 

 By Cauchy-Swharz inequality, the first integral can be majorated by $K|t'-t|^{1/2}\|a'\|_{L^2}$, the second integral is smaller than  $K (\rho^2(t-s)+\int_t^s|Y_u-Y'_u|^2 a_u du$). For the third integral write $c'_u(Y'_u)-c_u(Y_u)= b(u,Y_u)-b(u,Y'_u)-(\sigma(u,Y_u)-\sigma(u,Y'_u))\alpha'_u+ \sigma(u,Y_u)(\alpha_u-\alpha'_u)$ and we use ($\psi ii$), and eventually we use ($\psi iii$) in order to estimate the to  last integrals. Hence we have
\begin{eqnarray*}|Y_s- Y'_s|^2 &\leq & K(|x-x'|^2+(t-t')^{1/2}\|a'\|_{L^2}^2+\rho^2)+ K\int_t^s 
|Y_u-Y'_u|^2 a_udu\\
&&+ K\rho^2(s-t)^{1/2}\|a'\|_{L^2}+ K\int_t^s|Y_u-Y'_u|^2 a'_udu 
\\
&&+\left|\int_t^sD_{y}\psi_{1,\rho}(Y_u,Y'_u)\cdot\sigma(u,Y_u)  (\alpha_u-\alpha'_u)du\right|
\end{eqnarray*}
or equivalently
\[|Y_s- Y'_s|^2 \leq K g_\rho(s)+K\int_t^s|Y_u-Y'_u|^2 (1+|\alpha_u|+|\alpha'_u|)du\]
for some positive constants $C_1$ and $C_2$. By Gronwall's lemma the proof is complete. \quad $\diamond$

Fix $\alpha\in L^2$ and $x,x'\in \oGO$, and apply Lemma \ref{Li} to $Y=Y^{t,x,\alpha}$ and to $Y'=Y^{t',x',\alpha}$. 
We have, have $g_\rho(s)= \rho^2+|x-x'|^2+|t-t'|^{1/2}$, which gives 
$|Y_s-Y'_s|^2\le K(\rho^2+|x-x'|^2 +|t-t'|^{1/2})$. Letting $\rho$ to 0, we have obtained the desired result.

Fix $s_0\in [t,T]$  and $\alpha\in L^2(0,T;\R^d)$. Consider $Y:s\mapsto Y^{t,x,\alpha}_s$ and $Y': s\mapsto Y^{t,x,\alpha}_{s_0}$. 

\subsection*{Proof of (3)} 

Let us first prove that for each $(t,x)\in[0,T]\times \oGO$ and for each $s\in[t,T]$ one has
\begin{equation}\label{holder-x}
\sup_{u\in[s,t]}|Y^{t,x,\alpha}_u-x|\le K\sqrt{t-s}.
\end{equation}
Indeed, define  $Y$ and $Y'$ by  $Y_u=Y^{t,x,\alpha}_u$ and $Y'_u=x$ for $u\in[t,T]$. Put $c(u,y)=b(u,y)-\sigma(s,y)\alpha_s$ and $c'(u,y)=0$ for $(u,y)\in[t,T]\times \oGO$. The same computation as in Lemma \ref{Li}  gives,  for all $s'\in[t,s]$, 
 $$\psi_{1,\rho}(Y_{s'},x)\le \int_t^{s'} D_x\psi_{1,\rho}(Y_u,x)c(u,Y_u)du\le \int_t^sK |Y_u-x|(1+|\alpha_u|)du,$$ 
 which gives, using Cauchy-Shwarz inequality 
 $$|Y_{s'}-x|^2\le  K \left(\sup_{u\in[t,s]}  |Y_u-x|\right)\sqrt{t-s} \left(1+\|\alpha\|_{L^2(0,T)}\right).$$
 Passing to the supremum over $s'\in[t,s]$, we obtain
 $\displaystyle \sup_{u\in[t,s]}  |Y_u-x|\le K\sqrt{t-s}$.
 
 \vline
 
 Fix now $(t,x)\in[0,T]\times 
\oGO$ and $s,s'\in[t,T]$, with $s'\le s$. Put $x'=Y^{t,x}_{s'}$. By the previous result we have
$|Y^{s',x'}_s-x'|\le K\sqrt{s-s'}$, and by the uniqueness result we have the flow property $Y^{t,x}_s=Y^{s',x'}_s$, hence we have proved that $|Y^{t,x}_s-Y^{t,x}_{s'}|\le K \sqrt{s-s'}$.

\subsection*{Proof of (4)}
We apply Lemma \ref{Li} to $Y=Y^{t,x,\alpha}$ and $ Y^n=Y^{t,x,\alpha^n}$. Then $g^n_{\rho}$ is given  for all $s\in [t,T]$ by
$$g^n_{\rho}(s)=\rho^2+ \int_{t}^sf_\rho(u,Y_u,
Y^n_u)(\alpha_u-\alpha^{n}_u)du.$$
where $f_\rho(t,y,y')=D_x\psi_{1,\rho}(y,y')\sigma(t,y)$ is continuous on $ [0,T]\times\oGO$ $\times$ $\oGO$.\\
We first prove that a subsequence of $g^n_{\rho}$ converges pointwise to $\rho^2$ as
$n$ goes to $\infty$. We remark, by assertion (3) of Theorem \ref{leY}, that the sequence  $(Y^n)$ is bounded in
$\C^{0,1/2}([0,T];\oGO)$ and therefore is relatively compact.
Let $\overline Y$ be one of its limit in $X$. Let us prove that this limit is $Y$. Extracting a subsequence if necessary, one can suppose that the sequence $(Y^{n_p})_p$ converges to $\overline Y$. We
write
\begin{eqnarray*}g^{n_p}_{\rho}(s)&=&\rho^2+ \int_{t}^s(f_{\rho}(u,Y_u,Y^{n_p}_u)
-f_{\rho}(u,Y_u,\overline Y_u))(\alpha_u-\alpha^{n_p}_u)du\\
&&+ \int_{t}^s f_\rho(u,Y_u,\overline Y_u)\cdot
(\alpha_u-\alpha^{n_p}_u)du.
\end{eqnarray*}
The first integral converges to 0 as $p$ goes to $\infty$ by Lebesgue's Theorem  and the second integral
converges to $0$ by definition
of the weak convergence of $(\alpha^n)$ to $\alpha$. Now as $(\alpha^n)$ is bounded in $L^2(0,T;\R^d)$ there exists  $K>0$  (independent of $n$ and $\rho$) such that for all $n\in \N$ and for all $s$ $\in$ $[0,T]$,  $$\displaystyle{ g^{n}_\rho(s)\leq \rho^ 2+K}.$$
It follows, applying again Lebesgue's Theorem in the inequality given by Lemma \ref{Li},  that  for all $\rho\in(0,1)$ and for all $s\in$ $(t,T)$,
$$\lim_{p \to \infty}|Y_s-Y^{n_p}_s|^2\leq \rho^2(1+ \int_t^s |\alpha_u|du), $$ 
Letting $\rho$ to $0$, we deduce that $(Y^{n_p})_{p}$ converges pointwise, and
even uniformly to $Y$ and by uniqueness of the limit we have $Y= \overline Y$. This implies that the whole sequence $(Y^n)_{n}$ converges uniformly to $Y$ and the proof of (A3) is complete.

The proof of Theorem \ref{leY} is now complete. \quad $\diamond$

\section*{Appendix C: discontinuous mixed  single or multiple optimal stopping problems}

In this appendix we  first study a mixed optimal control-optimal stopping time problem and we prove that a particular value function is the maximal viscosity supersolution of  a  variational inequality. Then we  prove a reduction result: the value function of a mixed optimal control-optimal {\em multiple} stopping problem can be writen as the value function of a mixed optimal control-optimal {\em single} stopping problem  with a new reward defined recursively.

\subsection*{C.1. A deterministic mixed optimal control-optimal single stopping problem}
We first  study the following mixed optimal control-optimal single stopping problem.
For each bounded borelian real valued function $\psi$ defined on $[0,T]\times \oGO$  and for each $(t,x)$ in $[0,T]\times \oGO$  define the value function 
 $V[\psi](t,x)$ by  
\begin{equation}\label{vinfsup}
V[\psi](t,x)=\inf_{\alpha\in L^2(t,T)}\sup_{\theta\in[t,T]}\left\{
\frac{1}{2}\int_t^\theta
|\alpha_s|^2ds+\psi(\theta,Y^{t,x,\alpha}_\theta)\right\}, 
\end{equation}
where $Y^{t,x,\alpha}$ is the unique solution of (\ref{Y}).

When $\psi$ is upper-semicontinuous (usc), we show that  this value function is caracterized as the maximal viscosity subsolution of the following equation
\begin{equation}\label{Einfsup}\left\{\begin{array}{l}
{\displaystyle  \min\left(-\dt{V}+\frac{1}{2}|\sigma^TD V|^2-b\cdot
D V, V- \psi\right)=0\mbox{ in 
}[0,T)\times\GO}\vspace{2mm}\\
{\displaystyle \dpart{ V}{\gamma}=0 \mbox{ in }[0,T)\times 
\dGO,\
 V(T)= \psi(T)\mbox{ on }\oGO}
\end{array}\right.
\end{equation}
The proof follows different results of Barles and Perthame \cite{BP1}. We adapt them here to our context.

\begin{lemma}\label{L1}
$V[\psi^*]$ is usc and $V[\psi_*]$ is lsc. In particular, if $\psi$ is continuous, $V[\psi]$ is continuous.
\end{lemma} 
\noindent{\em Proof: } {\em Step 1:}
Suppose, by contradiction, that $V[\psi^*]$ is not usc. Then there exist $(t,x)\in [0,T]\times \oGO$, a sequence  $(t_n,x_n)_n$  that converges to $(t,x)$ and $\eps>0$ such that,
\begin{equation}\label{E1}
V[\psi^*](t,x) +2\eps \le \lim_{n\to \infty} V[\psi^*](t_n,x_n).
\end{equation}
Now from the one hand there exists $\overline \alpha$ such that for all  $ \theta\in[t,T],$ 
\begin{equation} \label{E2}
\frac{1}{2}\int_t^{\theta}
|\overline\alpha_s|^2ds+\psi^*(\theta,Y^{t,x,\overline\alpha}_{\theta})\le V[\psi^*](t,x)+\eps.
\end{equation}
From the other hand for each $n\in \N$ there exists $\theta_n\in[t_n,T]$ such that 
$$  V[\psi^*](t_n,x_n) \le \frac{1}{2}\int_{t_n}^{\theta_n}
|\overline\alpha_s|^2ds+\psi^*(\theta_n,Y^{t_n,x_n,\overline\alpha}_{\theta_n}).$$
Extracting a sequence if necessary, we have that $\theta_n$ converges to $\overline \theta \in[t,T]$.   
By the regularity of $Y$ given by Appendix B, we obtain that 
\begin{equation}\label{E3}\lim_{n\to \infty}  V[\psi^*](t_n,x_n) \le\frac{1}{2}\int_t^{\overline\theta}
|\overline\alpha_s|^2ds+\psi^*(\overline \theta,Y^{t,x,\overline\alpha}_{\overline\theta}).
\end{equation}
Now by (\ref{E2}) with $\theta= \overline \theta$ and by (\ref{E1}) and (\ref{E3}) we obtain
$$\frac{1}{2}\int_t^{\overline\theta}
|\overline\alpha_s|^2ds+\psi^*(\overline \theta,Y^{t,x,\overline\alpha}_{\overline\theta})+ \eps\le V[\psi^*](t,x) +2\eps \le \frac{1}{2}\int_t^{\overline\theta}
|\overline\alpha_s|^2ds+\psi^*(\overline \theta,Y^{t,x,\overline\alpha}_{\overline\theta}), $$
hence $\eps\le 0 $,  which is  the expected contradiction. 

\vline

{\em Step 2:} Suppose by contradiction that $V[\psi_*]$ is not lsc. Then there exist $(t,x)\in [0,T]\times \oGO$, a sequence  $(t_n,x_n)_n$  that converges to $(t,x)$ and $\eps>0$ such that,
\begin{equation}\label{E1b}
V[\psi_*](t,x)  \ge \lim_{n\to \infty} V[\psi_*](t_n,x_n)+ 3\eps.
\end{equation}
Note that for each $n$ there exists $\alpha^n\in L^2$ such that 
\begin{equation}\label{E2b}
V[\psi_*](t_n,x_n)+\eps\ge \sup_{\theta\in[t_n,T]}\frac{1}{2}\int_{t_n}^{\theta}
|\alpha^n_s|^2ds+\psi_*(\theta,Y^{t_n,x_n,\alpha^n}_{\theta})
\end{equation}
hence the sequence $(\alpha^n)$ is bounded in $L^2$ and extracting a subsequence if necessary, we can suppose that it converges to $\overline \alpha$ weakly in $L^2$. 

Now, as $\displaystyle \sup_{\theta\in[t,T]}\frac{1}{2}\int_{t}^{\theta}
|\overline\alpha_s|^2ds+\psi_*(\theta,Y^{t,x,\overline\alpha}_{\theta})\ge V[\psi_*](t,x)$ there exists $\overline \theta\in[t,T]$ such that
\begin{equation}\label{E3b}
\frac{1}{2}\int_t^{\overline\theta}
|\overline\alpha_s|^2ds+\psi_*(\overline\theta,Y^{t,x,\overline\alpha}_{\overline\theta})+\eps\ge V[\psi_*](t,x),
\end{equation}
For each $n\in\N$ define $\theta_n=t_n\wedge \overline \theta$. One has $\theta_n\in[t_n,T]$ and as $t_n\to t$ one has $\theta^n\to \overline\theta$. By (\ref{E2b}) we have
$\displaystyle V[\psi_*](t_n,x_n)+\eps\ge  \frac{1}{2}\int_{t_n}^{\theta_n}
|\alpha^n_s|^2ds+\psi_*(\theta_n,Y^{t_n,x_n,\alpha^n}_{\theta_n}).$
Adding $2\eps$ and passing to the liminf  we obtain, by continuity of $Y$ and lower semicontinuity of $\psi_*$
\begin{equation}\label{E4b} \lim_{n\to \infty} V[\psi_*](t_n,x_n)+3\eps\ge  \frac{1}{2}\int_{t}^{\overline \theta}
|\overline\alpha_s|^2ds+\psi_*(\overline\theta,Y^{t,x,\overline\alpha}_{\overline\theta})+2\eps .
\end{equation}
Now by (\ref{E3b}), (\ref{E1b}) and (\ref{E4b})  we obtain
$$ \frac{1}{2}\int_t^{\overline\theta}
|\overline\alpha_s|^2ds+\psi_*(\overline\theta,Y^{t,x,\overline\alpha}_{\overline\theta}) +\eps\ge  V[\psi_*](t,x) \ge \frac{1}{2}\int_t^{\overline\theta}
|\overline\alpha_s|^2ds+\psi_*(\overline\theta,Y^{t,x,\overline\alpha}_{\overline\theta}) +2\eps, $$
and the expected contradiction  $0\ge \eps$ follows.
\quad $\diamond$

\vline

\begin{lemma}\label{L2}
$V[\psi^*]$ (resp. $V[\psi_*]$) is a viscosity subsolution  (resp. supersolution) of  (\ref{Einfsup}) with obstacle $\psi$. In particular, if $\psi$ is continuous, $V[\psi]$ is a continuous solution of (\ref{Einfsup}) with obstacle $\psi$.
\end{lemma}

{\small 
For completness let us recall the definition of a viscosity subsolution and supersolution  of equation (\ref{Einfsup}). For simplicity we define $H(D\varphi)(t,x)= \frac{1}{2}|\sigma^TD\varphi(t,x)|^2-b\cdot D \varphi(t,x)$
\begin{definition}
An usc locally bounded function $v$ defined on $[0,T]\times \oGO$ is a viscosity subsolution of equation (\ref{Einfsup}) if and only if 
 
 $\forall \varphi \in \C^1([0,T]\times \oGO)$, if $(t_0,x_0)\in [0,T]\times \oGO$ is a local minimum of $v-\varphi$, then
\\
1.\, if $(t_0,x_0)\in [0,T]\times \GO$, \, $ \min \left(-\dt{\varphi}+H(D\varphi), v- \psi^*\right)(t_0,x_0)\le 0,$
 \\
2.\, if $(t_0,x_0)\in [0,T]\times\dGO$, \,   $\min \left(-\dt{\varphi}+H(D\varphi), v- \psi^*, \dpart{ \varphi}{\gamma}\right)(t_0,x_0) \le 0.$
\\
 
A lsc locally bounded function $u$ defined on $[0,T]\times \oGO$ is a viscosity supersolution of equation (\ref{Einfsup}) if and only if 
 
 $\forall \varphi \in \C^1([0,T]\times \oGO)$, if $(t_0,x_0)\in [0,T]\times \oGO$ is a local maximum of $u-\varphi$, then
\\
1.\, if $(t_0,x_0)\in [0,T) \times \GO$, \, $ \min \left(-\dt{\varphi}+H( D \varphi), u- \psi_*\right)(t_0,x_0)\ge 0,$
\\
2.\,  if $(t_0,x_0)\in \{T\} \times \GO$, \, $  v(t_0,x_0)- \psi_*(t_0,x_0) \ge 0,$\vspace{2mm}
\\ 
3.\,  if $(t_0,x_0)\in [0,T)\times\dGO$, $\max\left( \min(-\dt{\varphi}+H( D \varphi), u- \psi_*), \dpart{ \varphi}{\gamma}\right)(t_0,x_0)\ge 0,$
\\ 
4.\,    if  $(t_0,x_0)\in\{T\} \times \dGO $, \, $\max\left(  u(t_0,x_0)- \psi_*(t_0,x_0), \dpart{ \varphi}{\gamma}(t_0,x_0)\right) \ge 0.$ 
\end{definition}}
\noindent{\em Proof: } Let us first recall the Dynamic Programming Principle, which proof is well known in the  deterministic case, even for a discontinuous  reward. 

For each $(t,x)\in[0,T]\times \oGO$ and for each $\tau\in[t,T]$ we have
\begin{equation}\label{PPD}
V[\psi](t,x)= \inf_{\alpha\in L^2} \sup_{\theta\in[t,T]} \left\{ \frac{1}{2}\int_t^{\theta\wedge \tau}
|\alpha_s|^2ds+\psi( \theta,Y^{t,x,\alpha}_{\theta}) {\bf 1}_{\{\theta<\tau\}}\right.
\end{equation}
 $$\hspace{4cm}\left.+ V[\psi](\tau, Y^{t,x,\alpha}_\tau){\bf 1}_{\{\theta\ge\tau\}}\right\}.$$

\noindent{\em Step 1: } Let us first prove that $V[\psi^*]$ is a viscosity subsolution.

 Let $\varphi\in \C^1([0,T]\times \oGO)$ and suppose that $(t_0,x_0)$ is a local maximum of $V[\psi^*]-\varphi$. Without loss 
 of generality we can suppose that $V[\psi^*](t_0,x_0)=\varphi(t_0,x_0)$, and we fix $r>0$ such that for all 
 $(t,x)\in [0,T]\times\oGO$ if  $|t-t_0|<r$ and $|x-x_0|<r$  then  $V[\psi^*](t,x)\le \varphi(t,x)$.\\
 
If 
 $\left[V[\psi^*](t_0,x_0)\le \psi^ *(t_0,x_0)\mbox{ or }\left(x_0\in \dGO\; \mbox{ and }\frac{ \partial\varphi}{\partial\gamma}(t_0,x_0) \le 0\right)\right]$ there
 is nothing to prove. Suppose  that  
  \\\mbox{} \hspace{0,7cm} $\left[V[\psi^*](t_0,x_0)> \psi^*(t_0,x_0) \mbox{ and   } \left(\mbox{if }x_0\in \dGO\mbox{ then  }\frac{ \partial\varphi}{\partial\gamma}(t_0,x_0) > 0\right)\right].$\\
   Fix $\alpha\in L^2(t,T; \R^m)$ a constant control, and denote by $\overline\theta$ a real number in $[t_0, T]$ which maximizes $ \frac{1}{2} \int^\theta_{t_0}|\alpha|^2 du+\psi^*(\theta,Y^{t_0,x_0,\alpha}_\theta)$ on $[t_0,T]$. One has
$\varphi(t_0,x_0)=V[\psi^*](t_0,x_0)\le \frac{1}{2}(\overline\theta-t_0)|\alpha|^2 +\psi^*(\overline\theta,Y^{t,x,\alpha}_{\overline\theta}),$
and  in particular, $\overline\theta>t_0$. Consider also $\tilde \theta=\inf\{ s\in [t_0,T], |Y^{t_0,x_0,\alpha}_s-x_0|\ge r\}.$ One has $\tilde \theta>t_0$.  Fix $\tau\in(t_0,T] $ such that $\tau < \min(\overline\theta,\tilde \theta,  t_0+ r)$.
By the Dynamic Programming Principle   one has 
$\varphi(t_0,x_0) =V[\psi^*](t_0,x_0)\le \frac{1}{2}(\tau-t_0)|\alpha|^2+ V[\psi^*](\tau,Y^{t_0,x_0,\alpha}_\tau)\le \frac{1}{2}(\tau-t_0)|\alpha|^2+ \varphi(\tau,Y^{t_0,x_0,\alpha}_\tau) .$
Substracting $\varphi(t_0,x_0)$, dividing by $-h=t_0-\tau<0$, and letting $h$ to 0 we obtain, for each $\alpha\in \R^m$, 
$$- \frac{\partial \varphi}{\partial t}(t_0,x_0) - D\varphi(t_0,x_0)\cdot b(t_0,x_0) - \left(\frac{1}{2}|\alpha|^2-D\varphi(t_0,x_0)\cdot\sigma(t_0,x_0) \alpha\right)\le 0.$$
If we chose $\alpha= (-D\varphi\cdot\sigma)(t_0,x_0)$ we obtain 
$\left(- \frac{\partial \varphi}{\partial t} +H( D\varphi)\right)(t_0,x_0)\le 0,$
hence we have proved that $V[\psi^*]$ is a viscosity subsolution. 

\vline

\noindent{\em Step 2:} Let us prove now that $V[\psi_*]$ is a viscosity supersolution.

 Let $\varphi\in \C^1([0,T]\times \oGO)$ and suppose that $(t_0,x_0)$ is a local minimum of $V[\psi_*]-\varphi$. Without loss 
 of generality we can suppose that $V[\psi_*](t_0,x_0)=\varphi(t_0,x_0)$,
and fix $r>0$ such that for all 
 $(t,x)\in [0,T]\times\oGO$ if  $|t-t_0|+|x-x_0|<r$  then  $V[\psi^*](t,x)>\varphi(t,x)$.\\ 
  Clearly, if  $x_0\in \dGO$  and one has  $\frac{ \partial\varphi}{\partial\gamma}(t_0,x_0) \ge 0$ there 
 is nothing to prove. We suppose that,  if  $x_0\in \dGO$  one has $\frac{ \partial\varphi}{\partial\gamma}(t_0,x_0) < 0$. For each $\alpha\in L^2(t_0,T)$   one has clearly,
 $\sup_{\theta\in[t_0,T]}\frac{1}{2}\int_{t_0}^\theta|\alpha_s|^2ds+ \psi_*(\theta,Y^{t_0,x_0,\alpha}_\theta) \ge \psi_*(t_0,x_0).$ Taking the infimum over $\alpha\in L^2(t_0,T)$ we have 
 $V[\psi_*](t_0,x_0)\ge \psi_*(t_0,x_0).$

It remains to show that, if  $t_0\neq T$ then
\begin{equation}\label{L2E2}
-\frac{\partial \varphi}{\partial t}+H(D\varphi)
(t_0,x_0)\ge 0.\end{equation}
Fix $\tau\in (t_0,T]$. The Dynamic Programming Principle gives
\begin{eqnarray*}
V[\psi_*](t_0,x_0)&=&\inf_{\alpha\in L^2(t_0,T)} \max\left[\sup_{\theta\in[t_0,\tau)} \frac{1}{2}\int_{t_0}^ \theta |\alpha_s|^2ds+ \psi_*(\theta,Y^{t_0,x_0,\alpha}_\theta) \right. \\
&& \quad \quad \quad \quad \quad \quad \quad  +\left.\frac{1}{2}\int_{t_0}^ \tau |\alpha_s|^2ds+ V[\psi_*](\tau,Y^{t_0,x_0,\alpha}_\tau)  \right],
\end{eqnarray*}
hence
$$V[\psi^*](t_0,x_0)\geq \inf_{\arglim{\alpha\in L^2(t_0,T)}{ \|\alpha\|_{L^2}\leq K}}\left[   \frac{1}{2}\int_{t_0}^ \tau |\alpha_s|^2ds+ V[\psi_*](\tau,Y^{t_0,x_0,\alpha}_\tau)\right],$$
where $K=V[\psi^*](t_0,x_0)$. By the regularity of $Y$ with respect to $\alpha$ and $\tau$, there exists $h>0$ such that for all $\alpha\in L^2(t_0,T)$ with $\|\alpha\|_{L^2}\le K$, for all $\tau\in[t_0,t_0+h]$ one has $|\tau-t_0|+	|Y^{t_0,x_0,\alpha}_\tau-X_0|<r$.
Therefore, for each $\tau\in(t_0,t_0+h]$ we have 
\begin{equation}\label{L2E4}\varphi(t_0,x_0)=V[\psi^*](t_0,x_0)\geq \inf_{\alpha\in L^2(t_0,T)}\left[   \frac{1}{2}\int_{t_0}^ \tau |\alpha_s|^2ds+ \varphi(\tau,Y^{t_0,x_0,\alpha}_\tau)\right].
\end{equation}
 Suppose by contradiction that (\ref{L2E2}) is not satisfied. Then, there exists $\eps>0$ such that $\frac{\partial \varphi}{\partial t}
(t_0,x_0) +bD\varphi (t_0,x_0)-\frac{1}{2}|\sigma^TD\varphi|^2(t_0,x_0) \ge 2\eps$. Taking $r$ smaller if necessary, we can suppose that for all $(t,x)\in[0,T]\times \oGO$, if $|t-t_0|+|x-x_0|<r$ then
$$\frac{\partial \varphi}{\partial t}
(t,x) +bD\varphi (t,x)-\frac{1}{2}|\sigma^TD\varphi|^2(t,x) \ge \eps .$$
 Changing the value of $r$ if necessary, we can suppose  in the case when  $x_0\in\dGO$, that for all $(t,x)\in[0,T]\times \oGO$, if $|t-t_0|+|x-x_0|<r$ then $\frac{\partial \varphi}{\partial \gamma}(t,x)<0$ and in the case when  $x_0\in\GO$  that if $|x-x_0|<r$ then $x\in \GO$. In both cases,  if $\tau\in[t_0,t_0+h]$ then for all $\alpha\in L^2(t_0,T)$ with $\|\alpha\|_{L^2}\le K$ we have  ${\bf 1}_{\dGO}(Y^{t_0,x_0,\alpha}_\tau)\frac{\partial \varphi}{\partial \gamma}(\tau,Y^{t_0,x_0,\alpha}_\tau)\le 0$.

 Noticing that $\displaystyle \inf_{\alpha\in \R^m}\frac{1}{2}|\alpha|^2-D\varphi\sigma(t,x)\alpha= -\frac{1}{2}|\sigma^TD\varphi(t,x)|^2$, we have  
 for $\tau\in(t_0,t_0+h]$, \\
 $\displaystyle\varphi(\tau,Y^{t_0,x_0,\alpha}_\tau)\ge \varphi (t_0,x_0)+ \int_{t_0}^\tau (\frac{\partial \varphi}{\partial t} +D\varphi b-D\varphi \sigma \alpha)(u,Y^{t_0,x_0,\alpha}_u) du  \ge  \varphi (t_0,x_0)+ \eps (\tau-t_0)-\frac{1}{2}\int_{t_0}^\tau|\alpha_u|^2du.$ Now (\ref{L2E4})  gives
$$
\varphi(t_0,x_0)\ge  \varphi (t_0,x_0)+ \eps (\tau-t_0)>\varphi(t_0,x_0).$$
 which provides the expected contradiction. \quad $\diamond$

\begin{lemma}\label{L3}
Let $(\psi_n)$ be a nonincreasing sequence of continuous functions on $[0,T]\times \oGO$ such that 
$\psi^*= \lim \downarrow \psi_n$. Then $V[\psi^ *]= \lim \downarrow V[\psi_n]$.
\end{lemma}
\noindent{\em Proof: } Suppose that $\psi_n \downarrow \psi^*$.  Clearly, $V[\psi_n]\geq V[\psi^*]$ and the sequence $V[\psi_n]$ is nonincreasing hence we obtain $V[\psi^ *]\le  \lim \downarrow V[\psi_n]$. Let us prove the second inequality. Fix $(t,x)\in[0,T]\times \oGO$, and $\eps>0$. There exists $\alpha^*\in L^2$ and $\theta^*\in[t,T]$ such that for all $\theta\in[t,T]$ one has
\begin{equation}\label{EV}
V[\psi^*](t,x) + \eps \ge \frac{1}{2}\int_t^{\theta^*}
|\alpha_s^*|^2ds+\psi^*( \theta^*,Y^{t,x,\alpha^*}_{\theta^*})\ge \frac{1}{2}\int_t^{\theta}
|\alpha^*_s|^2ds+\psi^*( \theta,Y^{t,x,\alpha^*}_{\theta}).
\end{equation}
Now, for each $n\in \N$, there exists $\theta_n\in[t,T]$ such that
$ V[\psi_n](t,x)\le \frac{1}{2}\int_t^{\theta_n}
|\alpha^*_s|^2ds+\psi_n( \theta_n,Y^{t,x,\alpha^*}_{\theta_n}).$
Extracting a sequence if necessary, we can suppose that $\theta_n$ tends to $\overline \theta$.
Fix $p\in \N$. For each $n\geq p$, we have 
$$ V[\psi_n](t,x)\le \frac{1}{2}\int_t^{\theta_n}
|\alpha^*_s|^2ds+\psi_p( \theta_n,Y^{t,x,\alpha^*}_{\theta_n}). $$
Letting $n$ to $\infty$,  
$  \lim_{n\to \infty} V[\psi_n](t,x)\le \frac{1}{2}\int_t^{\overline\theta}
|\alpha^*_s|^2ds+\psi_p(\overline  \theta,Y^{t,x,\alpha^*}_{\overline \theta}).$
Now passing to the limit in $p$ and using (\ref{EV}) we obtain
$  \lim_{n\to \infty} V[\psi_n](t,x) \le V[\psi^*](t,x) + \eps.$\quad $\diamond$
\begin{proposition}\label{Puu}
$V[\psi^*]$ is the maximal usc viscosity subsolution  of (\ref{Einfsup}) with obstacle $\psi$.
\end{proposition}
\noindent{\em Proof: } In view of Lemmas \ref{L1} and \ref{L2}, the only point which is left to show is the maximality of the solution. Let $v$ be a usc function which is a viscosity  subsolution  of (\ref{Einfsup}) with obstacle $\psi$. Let $\psi_n$ be a nonincreasing sequence of continuous functions on $[0,T]\times \oGO$ such that 
$\psi^*= \lim \downarrow \psi_n$. 
Since $\psi^n\geq \psi^*$,  $v$ is also a viscosity subsolution of equation
(\ref{Einfsup}) with obstacle $\psi^n$.
By  Lemma \ref{L2},  $V[\psi^n]$ is a continuous   viscosity solution of the same equation, and  $ v\le V[\psi^n]$ by Theorem \ref{TU}. As by Lemma \ref{L3}, $V[\psi^ *]= \lim \downarrow V[\psi_n]$,  we obtain  $ v\le V[\psi^*]$.\quad $\diamond$

\subsection*{C.2. Reduction of  multiple stopping to  single stopping problems}
Let $(\psi_i)_{i\in \N}$ be a family of real valued bounded measurable functions defined on $[0,T]\times \oGO$ and consider for each nonempty finit subset $I$ of $\N$ and for each $(t,x)\in[0,T]\times \oGO$ the mixed optimal control--optimal {\em multiple} stopping problem 
\begin{equation}\label{v[I]} v^I(t,x)=\inf_{\alpha \in L^2[t,T]}\inf_{\theta_I\in
[t,T]^N}\left\{\int_t^{\vee_{i\in I}\theta_i}|\alpha_s|^2ds
+\sum_{i\in I}\psi_i(\theta_i,Y^{t,x,\alpha}_{\theta_i})\right\}.
\end{equation}
the value function of the following mixed optimal control--optimal {\em single} stopping problem
\begin{equation}\label{u[I]}
 u^{I}(t,x)=\inf_{\alpha\in L^2}\inf_{\theta\in[t,T]}\left\{
\frac{1}{2}\int_t^\theta
|\alpha_s|^2ds+\phi(\theta, Y^{t,x,\alpha}_\theta)\right\},
\end{equation}
where the {\em new reward} is defined recursively by
\begin{equation}\label{phiI}
\phi=\left\{\begin{array}{ll} \psi_i& \mbox{ if } I=\{i\} \mbox{ with } i\in \N,\\
\displaystyle \min_{i\in I}\left\{\psi_i+
v^{I\backslash\{i\},A}\right\}& \mbox{ if } I \mbox{ contains 2 or more elements.}
\end{array}\right.
\end{equation}
\begin{proposition}\label{Preduction} For each finit nonempty subset $I$ of $\N$  let $v^I$ be defined by (\ref{v[I]}) and let $u^ I$ be defined by (\ref{u[I]}) and (\ref{phiI}). Then for each $(t,x)\in[0,T]\times \oGO$
 one has $v^I(t,x)=u^I(t,x).$
\end{proposition}
\noindent{\em Proof: } Fix a nonempty finit subset $I$ of $\N$ of cardinal $N$.
When $I$ contains only one element, there is nothing to prove. Suppose now that $I$ has two or more elments. Let us prove first that for each $(t,x)$ $\in $ $[0,T]\times \oGO$, 
$v^I(t,x)\leq u^ I(t,x)$.\\
Fix $(t,x)\in [0,T]\times \oGO$. Consider a partition $(A_j)_{j\in I}$ of $[t,T]^N$ such that for each  $\theta_I\in[t,T]^N$ if $\theta_I\in A_j$ then one has $\theta_j= \wedge_{i\in I} \theta_i$.\\

Fix $\alpha \in L^2(t,T)$ and $\theta_I\in[t,T]^N$,
\begin{eqnarray*}\int_t^{\vee_{i\in I}\theta_i}|\alpha_s|^2ds
&\hspace{-3mm}+&\hspace{-3mm}\sum_{i\in I}\psi_i(\theta_i,Y^{t,x,\alpha}_{\theta_i})= \displaystyle \sum_{j\in I} {\bf 1}_{A_j}(\theta_I)
\left\{ \int_t^{ \theta_j} |\alpha_s|^2ds + \psi_j (\theta_j,Y^{t,x,\alpha}_{\theta_j})  \right. \\
&& +  \left.\int_{\theta_j}^{\vee_{i\in I\backslash \{j\}}\theta_i}|\alpha_s|^2ds
+\sum_{i\in I\backslash \{j\}}\psi_i(\theta_i,Y^{t,x,\alpha}_{\theta_i})
\right\} =(I).
\end{eqnarray*}
Clearly, by uniqueness for equation (\ref{Y}), the second terme of the right-hand side can be minorated by  $ v^{I\backslash \{j\}}(\theta_j,Y^{t,x,\alpha}_{\theta_j})$, hence
\begin{eqnarray*}(I)&\ge & \sum_{j\in J}{\bf 1}_{{A}_j}(\theta_I) \left\{ \int_t^{ \theta_j} |\alpha_s|^2ds + \psi_j (\theta_j,Y^{t,x,\alpha}_{\theta_j}) + v^{I\backslash \{j\}}(\theta_j,Y^{t,x,\alpha}_{\theta_j})\right\} \\
&\ge & \sum_{j\in J}{\bf 1}_{{A}_j}(\theta_I) \left\{ \int_t^{ \theta_j} |\alpha_s|^2ds + \phi(\theta_j, Y^{t,x,\alpha}_{\theta_j})\right\} \ge \left(\sum_{j\in J}{\bf 1}_{{A}_j}(\theta_I)\right) u^I(t,x).
\end{eqnarray*}
Hence $(I)\ge   u^I(t,x)$. Taking the infimum over  $\alpha \in L^2(t,T)$ and $\theta_I\in[t,T]^N$, we obtain $v^I(t,x)\geq u^I(t,x)$.\\

Let us now prove the reverse inequality. For simplicity, suppose first that there exist an optimal time $\theta^*\in[t,T]$ and an optimal 
control $\alpha^*\in L^2(t,T)$   for $u^I(t,x)$, and for each $i\in I$ there exist $\tilde \alpha^ {i*}\in L^2(\theta^*,T)$ and $\tilde \theta^{i*}  \in [\theta^ *,T]^{N-1}$ that are optimal for $v^{I\backslash\{i\}}(\theta^*, Y^{t,x,\alpha^*}_{\theta^*})$. Define $ \alpha^{i *}= \alpha^* {\bf 1}_{[t,\theta^*)}+ \tilde\alpha^{i*}{\bf 1}_{[\theta^*,T)}$. Note that $\alpha^{i*}$ is optimal for $u^I(t,x)$, that  $Y^{t,x,\alpha^*}_{\theta^*}=  Y^{t,x,\alpha^{i*}}_{\theta^*}$ and that  $\alpha^{i*}$ is also optimal for $v^{I\backslash\{i\}}(\theta^*, Y^{t,x,\alpha^*}_{\theta^*})$.

 Let $(B_i)_{i\in I}$ be a partition of $[t,T]$ such that  for each $s\in[t,T]$, if $s\in B_i$ then $\phi(s,Y^{t,x,\alpha^*}_s)= \psi_i(s,Y^{t,x,\alpha^*}_s)+ v^{I\backslash\{i\}}(s,Y^{t,x,\alpha^*}_s)$. 
\begin{eqnarray*} u^I(t,x)&=&\int_t^{\theta^*}|\alpha^*_s|^2ds + \phi(\theta^*, Y^{t,x,\alpha^*}_{\theta^*})\\
&=&  \sum_{i\in I}{\bf 1}_{B_i} (\theta^*)\left(\int_t^{\theta^*}|\alpha^*_s|^2ds +  \left\{\psi_i +  v^{I\backslash\{i\}}\right\}(\theta^*, Y^{t,x,\alpha^*}_{\theta^*})\right)\\
&=&  \sum_{i\in I}\left({\bf 1}_{B_i}(\theta^*)\left\{ \int_t^{\theta^*}|\alpha^{i*}_s|^2ds  +\psi_i(\theta^*,Y^{t,x,\alpha^{i*}}_{\theta^*})\right.\right.\\
&& \quad\quad \quad \left.\left. + \quad \int_{\theta^*}^{\vee_{j \in I\backslash\{i\}}\tilde \theta^{i*}_j}|\alpha^{i*}_s|^2ds  + \displaystyle \sum_{j\in I\backslash\{i\}} \psi_j(\tilde\theta^{i*}_j,Y^{t,x,\alpha^{i*}}_{\tilde\theta^{i*}_j})   \right\}\right) 
\end{eqnarray*}
For each $i$ define $ \theta^{i*}\in[t,T]^ N$ by  $\theta^{i*}_i=\theta^*_i$ and $ \theta^{i*}_j=\tilde \theta^{i*}_j$ for $j\in I\backslash\{i\}$.
One has now
$$ u^I(t,x)=\sum_{i\in I}{\bf 1}_{B_i}(\theta^*)\left\{  \int_t^{\vee_{j\in I} \theta^{i*}_j}| \alpha^{i*}_s|^2ds+ \sum_{j\in I}
\psi_j(\theta^{i*}_j,Y^{t,x,\alpha^{i*}}_{\theta^{i*}_j})\right\}
\ge   v^I(t,x).
$$
In general, there is no optimal stopping control and stopping times, but for each $\eps>0$ on can find $\eps/2$ optimal  $ \theta^*$ and $\alpha^*$ for $u^I(t,x)$ and 
$\eps/2$ optimal $\tilde \theta^{i*}_p$ and $\tilde\alpha^{i*}_p$ for $v^{I\backslash\{i\}}(\theta^*, Y^{t,x,\alpha^*}_{\theta^*})$. Building $\alpha^{i*}$ and $\theta^{*i}$ as above we obtain that for each $\eps>0$,
$u^I(t,x) +\eps \ge v^I(t,x).$ \quad $\diamond$

\section*{Appendix D: a strong comparison result}
In this appendix we prove a strong comparison result  for viscosity solutions of a first order  variational inequality with Neumann boundary conditions and with continuous obstacle.
\begin{theorem}
\label{TU}
Assume (\ref{gamma}) and (\ref{bsigma1})-(\ref{bsigma2}) and let
$\psi\in \C([0,T]\times \oGO; \R)$.
If $u, \, v:[0,T]\times\oGO\to\R$  are respectively
usc viscosity subsolution and lsc viscosity supersolution of
\begin{equation}\label{U D1}
\left\{\begin{array}{l}
{\displaystyle  \min\left(-\dt{w}+\frac{1}{2}|\sigma^TDw|^2-b\cdot
Dw,w-\psi\right)=0\mbox{ in }[0,T)\times\GO,}\vspace{2mm}\\
{\displaystyle \dpart{w}{\gamma}=0 \mbox{ in }[0,T)\times \dGO,\
w(T)=\psi(T)\mbox{ on }\oGO,}
\end{array}\right.
\end{equation}
then $u\leq v$ on $[0,T]\times\oGO$.
\end{theorem}

\vline

Note that the difficulty of proving this strong comparison result is double. First, we have to handle the Neuman condition, and the test function of Appendix A was built to that aim. Second, even though the equation is of first order and no Ishii lemma is needed, the quadratic term $|\sigma^T Dv|^2 $ has to be taken with care. 

\vline 

\noindent{\em Proof: } In the following we denote for all $\eps>0$ by $\Psi_\eps$ the function defined  by $\Psi_\eps=\psi_{\eps,\eps^2}$ where $\psi_{\eps,\rho}$ is the test-function of Lemma \ref{LT}.

\vline 

 Let $u$, $v$ be
respectively a bounded usc viscosity subsolution and a bounded lsc 
viscosity supersolution
of (\ref{U D1}). We first remark that the terminal condition holds in a stronger sense
that the viscosity sense.

\begin{proposition} \label{bord}
    For all $x\in\oGO$, $u(T,x)\leq \psi(T,x)$.
\end{proposition}

\noindent{\em Proof: } Fix $x_{0}\in\oGO$. For all $\eps>0$, put $\varphi_\eps(t,x)=\Psi_{\eps}(x,x_{0})+\frac{T-t}{\eps^2}$ and let
$(t_{\eps},x_{\eps})$ be a global
maximum of the usc function $ u-\varphi_{\eps}$.
Then
$u(T,x_{0})-\varphi_{\eps}(x_{0},x_{0})\leq
u(t_{\eps},x_{\eps})-\varphi_{\eps}(x_{\eps},x_{0})$
which implies, by Lemma \ref{LT},
$\frac{1}{2}\frac{|x_{\eps}-x_{0}|^2}{\eps^2}- 
K\eps^2+\frac{T-t_{\eps}}{\eps^2}
\le \varphi(t_\eps,x_\eps)
\le u(t_{\eps},x_{\eps})-u(T,x_{0})+\varphi_{\eps}(T,x_0)
\leq 2\|u\|+ K\eps^2.$
We deduce that $t_{\eps}$ and $x_{\eps}$ go respectively to $T$ and
$x_{0}$ as $\eps$ goes to 0 and, by the upper semicontinuity of $u$,
that, 
\begin{equation}\label{caconv}
\displaystyle{\lim_{\eps\to 0}u(t_{\eps},x_{\eps})= u(T,x_{0})}
\mbox{ and }
\displaystyle{\lim_{\eps\to 0 }\frac{|x_{\eps}-x_0|^2}{\eps^2}= 0}, \; 
\displaystyle{\lim_{\eps\to 0 }\frac{T-t_{\eps}}{\eps^2}= 0}.
\end{equation}
As, by Lemma \ref{LT}, $\frac{\partial \varphi_\eps}{\partial \gamma} (t,x) >0$ if $x\in\dGO$, by definition of viscosity subsolution, one has,  for all $\eps>0$ and
$x_{\eps}\in\oGO$,
$$ \min\left(\frac{1}{\eps^2}
+\frac{1}{2}|\sigma^TD_{x}\varphi_\eps|^2-b\cdot
D_{x}\varphi_\eps, u-\psi\right)(t_{\eps},x_{\eps})\leq 0.$$
Now, assume that, for some subsequence,
$u(t_{\eps},x_{\eps})>\psi(t_{\eps},x_{\eps})$. Then, necessarily, for
this subsequence,
\begin{equation}\label{impossible}
1\leq
\eps^2|D_{x}\Psi_{\eps}(x_{\eps},x_{0})|
(\frac{\|\sigma\|^2}{2}|D_{x}\Psi_{\eps}(x_{\eps},x_{0})|+\|b\|).
\end{equation}
But, by lemma \ref{LT},
$\eps|D_{x}\Psi_{\eps}(x_{\eps}-x_{0})|\leq 
K(\frac{|x_\eps -x_0|}{\eps}+\eps^3)$
and therefore (\ref{impossible}) cannot occur because of
(\ref{caconv}) and for all $\eps>0$,
$u(t_{\eps},x_{\eps})\leq\psi(t_{\eps},x_{\eps})$. We conclude by
letting $\eps$ go to 0 and using (\ref{caconv}).\quad $\diamond$

\vline

For all $0<\nu<1$ and all $\delta>0$, let
\[ M_{\nu,\delta}=\sup_{(t,x)\in[0,T]\times\oGO}(\nu
u(t,x)-v(t,x)-\delta (T-t)).\]
Our aim is to prove that $M_{\nu,\delta}\leq (1-\nu)\|\psi\|$ 
which will give the
conclusion of the theorem by letting $\nu$ and $\delta$ tend to 1 and
0 respectively. To do so, we define, for all $\eps$, $\alpha>0$,
$M^{\eps,\alpha}_{\nu,\delta}$ as being the supremum over
$[0,T]^2\times\oGO^2$ of the function
\[(t,s,x,y)\mapsto \nu u(t,x)-v(s,y)-\delta
(T-s)-\Psi_{\eps}(x,y)-\frac{(t-s)^2}{\alpha^2}\]
and denote by $(\ct,\cs,\cx,\cy)$ an optimal point (recall that $u$ 
and $v$
are bounded and respectively usc and lsc).

We first notice that, since $\Psi_{\eps}(x,x)=0$ for all $x\in\dGO$,
\[\nu u(T,x)-v(T,x)\leq M^{\eps,\alpha}_{\nu,\delta} \]
for all $x\in\dGO$, which implies by Lemma \ref{LT}
\begin{equation}\label{alpha}
\frac{1}{2}\frac{|\cx-\cy|^2}{\eps^2}\leq 
2(\|u\|+\|v\|)+ K
\eps^2\mbox{ and }\frac{(\ct-\cs)^2}{\alpha^2}\leq
2(\|u\|+\|v\|)+ K\eps^2.
\end{equation}
Therefore, up to some subsequences, $\ct$, $\cs$ and $\cx$, $\cy$ 
converge
respectively to some $\overline t$ and $\overline x$ in $[0,T]\times 
\oGO$ as $\alpha$
and $\eps$ go to 0.

Now we proceed as in Propostion \ref{bord}. For all $\alpha$, $\eps>0$
\[\nu u(\overline t,\overline x)-v(\overline t,\overline x)-\delta 
(T-\overline
t) -\Psi_{\eps}(\overline x,\overline x)\leq M^{\eps,\alpha}_{\nu,\delta},\]
so that, by Lemma \ref{LT}, the upper semicontinuity of
$u$  and the lower semicontinuity of $v$, as $\alpha$ and $\eps$ go 
to 0,
\begin{equation}\label{espace}
\frac{	|\cx-\cy|^2}{\eps^2}+\frac{(\ct-\cs)^2}{\alpha^2}\to 0.
\end{equation}
As a consequence, we get, as $\alpha$ and $\eps$ tend to 0,
\begin{equation}\label{total}
    M^{\eps,\alpha}_{\nu,\delta}\to M_{\nu,\delta}.
\end{equation}

We define, for all $(t,x),\ (s,y)\in [0,T]\times\oGO$,
\begin{eqnarray*}
    \varphi_1(t,x)&=&\myfrac{1}{\nu}
    \left( v(\cs,\cy)+\delta(T- \cs)+\Psi_\eps(x,\cy)+
    \myfrac{|t-\cs|^2}{\alpha^2}\right)\\
    \varphi_2(s,y)&=&
    u(\ct,\cx)-\delta(T- s)-\Psi_\eps(\cx,y)- 
\myfrac{|\ct-s|^2}{\alpha^2}
\end{eqnarray*}
and we apply the definition of viscosity solutions to $u$ and $v$:
$u-\varphi_1$ reaches its maximum at $(\ct,\cx)$ and when
$\cx\in\dGO$ we can check easily that
$ D\varphi_{1}(\ct,\cx)\cdot \gamma(\cx)>0$ by Lemma \ref{LT} and 
therefore the Neumann
boundary condition never holds. This imply that for all $\alpha$,
$\eps$,
\begin{equation}\label{pouru}
\min\left(-\dt{\varphi_{1}}+\frac{1}{2}
|\sigma^TD\varphi_1|^2-b\cdot\varphi_{1},
u-\psi\right)(\ct,\cx)\leq 0.
\end{equation}
For $v$, the situation is slightly different. As in the former case,
we deduce from Lemma \ref{LT} that the Neumann boundary condition
cannot hold when $\cy\in\dGO$, but if for some subsequence of
$(\alpha,\eps)$, $\cs=T$ then we can have $v(\cs,\cy)\geq
\psi(\cs,\cy)$ and no information on the partial differential
inequation. In this case, we remark that $\ct$ goes to $T$ (hence 
$\overline t$=T) and that, 
by
Proposition \ref{bord} and the upper semicontinuity of $u$, for all 
$\delta_{0}$,
$u(\ct,\cx)\leq u(T,\overline x) +\delta_0/2\leq \psi(T,\overline x) 
+\delta_0/2\leq \psi(\ct,\cx)+\delta_{0}$ for $\alpha$ and $\eps$
small enough. We deduce, from those two inequalities, by passing to 
the
limit as $\alpha$ and $\eps$ go to 0 and using (\ref{total}), that
\[M_{\nu,\delta}=\nu u(T,\overline x)-v(T,\overline
x)\leq \nu \psi(T,\overline x)-\psi(T,\overline x)+\delta_0\leq (1-\nu)\|\psi\|+\delta_{0}\nu\]
for all $\delta_{0}>0$, so that finally $M_{\nu,\delta}\leq 
(1-\nu)\|\psi\|$.

Now we are left with the case, when $\cs<T$ at least along a 
subsequence of $(\varepsilon,\alpha)$. We have
\begin{equation}\label{pourv}
\min\left(-\dt{\varphi_{2}}+\frac{1}{2}
|\sigma^TD\varphi_2|^2-b\cdot\varphi_{2},
v-\psi\right)(\cs,\cy)\geq 0.
\end{equation}

If, for some subsequence, $u(\ct,\cx)>\psi(\ct,\cx)$ then
(\ref{pouru}) and (\ref{pourv}) give respectively
\[-\dt{\varphi_{1}}(\ct,\cx)+\frac{1}{2}
|\sigma(\ct,\cx)^TD\varphi_1(\ct,\cx)|^2-b(\ct,\cx)\cdot\varphi_{1}(\ct,\cx)\leq 
0\]
and
\[-\dt{\varphi_{2}}(\cs,\cy)+\frac{1}{2}
|\sigma(\cs,\cy)^TD\varphi_2(\cs,\cy)|^2-b(\cs,\cy)\cdot\varphi_{2}(\cs,\cy)\geq 
0.\]
We multiply the first inequality by $\nu$ and substract the second 
one ;
we obtain a rather complicated inequality which has three kinds of 
terms:
the time derivative term,  the linear term and the quadratic term.

\noindent {\em The time derivative term } is the simplest one
\[-\nu \dt{\varphi_{1}}(\ct,\cx)+\dt{\varphi_{2}}(\cs,\cy)= \delta.\]
{\em The linear term } can be writen
$$(b(\cs,\cx)-b(\ct,\cx))\cdot
D_{x}\Psi_{\eps}(\cx,\cy)+(b(\cs,\cy)-b(\cs,\cx))\cdot
D_{x}\Psi_{\eps}(\cx,\cy) $$ $$ -b(\cs,\cy)\cdot
(D_{x}\Psi_{\eps}(\cx,\cy)+D_{y}\Psi_{\eps}(\cx,\cy))$$
and can be estimated, if $\omega_{b}$  and $K_{b}$ denote
respectively the modulus of
continuity with respect to $t$ and the Lipschitz constant with respect 
to $x$ of $b$ on $[0,T]\times\oGO$ and
by using Lemma \ref{LT}, by
\[ K\left( (\omega_{b}(|\ct-\cs|)+K_{b}|\cx-\cy|)\left(
\frac{|\cx-\cy|}{\eps^2}+\eps^2\right) +
\|b\|\left(\frac{|\cx-\cy|^2}{\eps^2}+\eps^2\right) \right).\]
We know, by (\ref{alpha}), that $|\ct-\cs|\leq C\alpha$ for some
constant $C$ independent of $\alpha$ and $\eps<1$, therefore, if we
choose $1>\eps>\omega_{b}(C\alpha)$ as $\alpha$ and $\eps$ go to 0,
this linear term goes to 0 by (\ref{espace}).

\noindent As far as {\em the quadratic term} is concerned, we first
remark that for all $ a$, $ b$ in $\R^m$ and all 
$0<\nu<1$,
\[\frac{1}{\nu}| a|^2-| b|^2\geq
-\frac{1}{1-\nu}| a+ b|^2\]
so that we are reduced to estimate
$$|\sigma(\ct,\cx)^TD_{x}\Psi_{\eps}(\cx,\cy)+\sigma(\cs,\cy)^TD_{y}\Psi_{\eps}(\cx,\cy)|^2,$$
which we do as for the linear term, concluding that it goes to 0 as
$\alpha$ and $\eps$ go to 0, providing that 
$1>\eps>\omega_{\sigma}(C\alpha)$.

In conclusion to all those estimates we obtain the contradiction
$\delta\leq 0$, and finally we necessarily have, for all $\alpha$
and $\eps>\omega_{b}(C\alpha)$ small enough, $u(\ct,\cx)\leq
\psi(\ct,\cx)$. This, combined with (\ref{pourv}) and (\ref{total}), 
yields
$M_{\nu,\delta}\leq (1-\nu)\|\psi\|$
and the proof is complete.\quad $\diamond$

\section*{Acknowledgements}
The author thanks Guy Barles for his relevant remarks and suggestions.

Address of Magdalena Kobylanski\\CNRS - UMR 8050 (LAMA)\\ Universit\'e de Marne-la-Vall\'ee\\5, boulevard Descartes\\Cit\'e Descartes - Champs-sur-Marne\\77454 Marne-la-Vall\'ee cedex 2 (FRANCE)\\magdalena.kobylanski@univ-mlv.fr
\end{document}